\documentclass[final]{siamltex}

\usepackage{amsmath}
\usepackage{mathrsfs}
\usepackage{amsfonts}
\usepackage{amssymb}
\usepackage{graphicx}
\usepackage{psfrag}
\usepackage{placeins}
\usepackage{epsfig}
\usepackage{amsopn}
\usepackage{amstext}
\usepackage{amssymb}
\usepackage{amscd}
\usepackage{latexsym}
\allowdisplaybreaks
\usepackage{geometry}
\usepackage{multirow}
\usepackage{diagbox}
\usepackage{color}
\usepackage{ulem}

\newtheorem{rem}[theorem]{Remark}

\newcommand{\bbm}{\begin{bmatrix}}
	\newcommand{\ebm}{\end{bmatrix}}

\renewcommand{\theequation}{\arabic{section}.\arabic{equation}}
\usepackage{subfigure}
\usepackage{graphicx}
\usepackage{float}
\usepackage[numbers,sort&compress]{natbib}

\begin{document}

\title{Parallel Domain Decomposition method for the fully-mixed Stokes-dual-permeability fluid flow model with Beavers-Joseph interface conditions}
\author{
	Zheng Li \thanks {%
		Department of Mathematics, Shanghai Maritime University, Shanghai, P.R.China. \texttt{ilizheng@msn.com}. }
\and Feng Shi \thanks{%
	College of Science, Harbin Institute of Technology, Shenzhen, P.R.China. \texttt{shi.feng@hit.edu.cn}. Partially supported by Foundation Research Project of Shenzhen (Grant No. GXWD20201230155427003-20200822102539001).}
\and Yizhong Sun \thanks{%
School of Mathematical Sciences, East China Normal University,
Shanghai, P.R.China. \texttt{bill950204@126.com}. }
\and Haibiao Zheng \thanks{%
School of Mathematical Sciences, East China Normal University,
Shanghai Key Laboratory of Pure Mathematics and Mathematical
Practice, Shanghai, P.R. China. \texttt{hbzheng@math.ecnu.edu.cn}.
Partially supported by   NSF of China (Grant No. 11971174) and NSF
of Shanghai (Grant No. 19ZR1414300) and Science and Technology
Commission of Shanghai Municipality (Grant No. 18dz2271000 and No. 19JC1420102). }}

 \maketitle

\begin{abstract}
In this paper, a parallel domain decomposition method is proposed for solving the fully-mixed Stokes-dual-permeability fluid flow model with  Beavers-Joseph (BJ) interface conditions. 
Three Robin-type boundary conditions and a modified weak formulation are constructed to completely decouple the original problem, not only for the free flow and dual-permeability regions but also for the matrix and microfractures in the dual-porosity media. We derive the equivalence between the original problem and the decoupled systems with some suitable compatibility conditions, and also demonstrate the equivalence of two weak formulations in different Sobolev spaces. Based on the completely decoupled modified weak formulation, the convergence of the iterative parallel algorithm is proved rigorously. To carry out the convergence analysis of our proposed algorithm,  we propose an important but  general convergence lemma for the steady-state problems. Furthermore, with some suitable choice of parameters, the new algorithm is proved to achieve the geometric convergence rate. Finally, several numerical experiments are presented to illustrate and validate the performance and exclusive features of our proposed algorithm.
\end{abstract}


\section{Introduction}
Multi-domain, multi-physics coupled problems are significant in many natural and industrial applications, such as the groundwater fluid flow in the karst aquifer,  petroleum extraction, industrial filtration, blood flow motion in the arteries, and so on. Up to now, a great deal of mathematical and physical models for free flow coupled with complicated porous media were constructed, including the Stokes-Darcy model \cite{DMQ02, Badea, CGHW08, Cao10}, dual-porosity-Stokes model \cite{HouJY16, Shan19, Mahbub2019, Mahbub2020, WangYS21, HouJY21, Shan22}, etc. Typically, we usually utilize the Darcy equation in the traditional Stokes-Darcy fluid flow system to simulate the single porosity model. Plenty of numerical methods for solving the Stokes-Darcy system can be referred to  Lagrange multiplier methods \cite{LMLayton2003, LMGatica03, LMGatica2012}, domain decomposition methods \cite{Discacciati07, Chen, Boubendir, Cao11, He15, Vassilev1, Sun2021, Sun2021tg}, optimized Schwarz methods \cite{Discacciati, Discacciati18, Gander}, and partitioned time stepping methods \cite{ Shan13, Layton20}, to name just a few. 

The Darcy equation in Stokes-Darcy system has some limitations in describing fluid flow in porous media with complicated geometrical structures. For example, a naturally fractured reservoir containing the multi-porosity/permeability regions is comprised of low permeable rock matrix blocks surrounded by an irregular network of natural microfractures \cite{Arbogast92, Douglas98, HouJY16, Shan19}. In order to describe the coupled flow in the dual-porosity media and conduits, Hou et al. \cite{HouJY16} constructed a new dual-porosity-Stokes model, and proposed four physically valid interface conditions to couple the Stokes equation and dual-porosity models on the interface, including a no-exchange condition, a mass balance condition, a force balance condition, and the Beavers–Joseph (BJ) conditions. Their  work 
well established the continuous model of the coupled dual-porosity-Stokes system to simulate multistage hydraulic fractured horizontal wellbore.

Obviously, one important task is how to develop effective numerical algorithms for the dual-porosity-Stokes equations. Inspired by the decoupled ideas for the Stokes–Darcy model, some researchers extended to study the Stokes-dual-permeability model. A natural decoupled method is the domain decomposition method (DDM), since it  can decouple the multi-domain, multi-physics problems naturally under the introduced interface boundary conditions \cite{Discacciati07, Chen, Boubendir, Cao11, He15, Vassilev1, Sun2021, Sun2021tg}, and there are many well studied off-the-shelf and efficient solvers for each decoupled subproblem. Based on the characteristics of easy-to-operate, high precision and convenient parallel computing, DDM has received extensive attention and applications undoubtedly. In \cite{HouJY21}, the authors extended the Robin-Robin 
DDM of the Stokes-Darcy model in \cite{Chen, Cao11} to deal with the dual-porosity-Stokes model. Noting that such methods were usually based on Galerkin approximation to two elliptic pressure equations in the dual-porosity media. Hou et al. \cite{HouJY21} studied the dual-porosity-conduit system with the Beavers-Joseph-Saffman-Jones (BJS) interface boundary conditions, and decoupled the Stokes equations and the dual-porosity model by DDM, but not for matrix systems and microfracture systems. They demonstrated the convergence of the proposed algorithms and got a geometric convergence rate with some suitable selections of the Robin parameters, and also utilized the optimized Schwarz methods to get the optimized Robin parameters, which can improve the convergence of the proposed algorithms.

In the study of multi-domain, multi-physics fluid flow, it is usually necessary not only to get accurate pressure results, but also to obtain accurate velocity information. 
Moreover, mixed finite element methods are usually superior to the classical Galerkin methods in many areas, due to the natures of discontinuity of the gradient of the pressure and continuity of the dual-permeability velocity in applications. To this end, in this paper, we try to design a Robin-type DDM for the fully-mixed Stokes-dual-permeability coupled model with more physically realistic BJ interface conditions. In order to completely decouple this steady-state model, the designed algorithm inevitably 
contains more explicit terms, which increase the difficulty of convergence analysis. We present an important but general convergence lemma for steady-state problems, which can overcome the difficulty of convergence analysis. To the best of our knowledge, this is the first work to completely decouple the fully-mixed Stokes-dual-permeability model without introducing any stabilization terms. By this mean, we can easily use many existing well-developed solvers or codes in a flexible way to solve two single dual-permeability equations and a single Stokes equation in parallel. Furthermore, we demonstrate the convergence of the parallel Robin-type DDM. Interestingly, our algorithm also has a mesh-independent geometric convergence rate with suitable choice of the Robin parameters. Noting that our proof begins with a rigorous equivalence analysis between the original and the modified weak formulations, so that the subsequent arguments can be only based on the completely decoupled modified weak formulation.

The rest of this paper is organized as follows. In Section 2, the fully-mixed Stokes-dual-permeability model with BJ interface conditions is described. In order to apply the idea of DDM to solve the fully-mixed Stokes-dual-permeability model, three Robin-type interface conditions are constructed and their equivalences with the original
interface conditions are proved under suitable  compatibility conditions. In Section 3, we introduce a modified weak formulation and rigorously demonstrate  its equivalence with the original weak formulation. Then a parallel Robin-type DDM and its convergence analysis are presented in Section 4. More importantly, the parallel DDM for the continuous model is proved to achieve geometric convergence by  appropriate choice of the Robin parameters. Finally three numerical tests are carried out to illustrate the exclusive features of our proposed DDM in Section 5. 

\section{Fully-mixed Stokes-dual-permeability Model with BJ Interface Conditions}
Consider two bounded domains denoted by $\Omega _{S},\Omega_{D}\subset R^{d}(d=2~\mathrm{or}~3)$ with an interface $\Gamma$ between two subdomains. Assume that the two bounded domains are non-overlapping, i.e., $\Omega _{S}\cap \Omega _{D}={\emptyset }$, and $\overline{\Omega }_{S}\cap\overline{\Omega }_{D}=\Gamma $. Denote $\overline{\Omega }=\overline{\Omega }_{S}\cup \overline{\Omega }_{D}$. Define  $\mathbf{n}_{S}$ and $\mathbf{n}_{D}$ as the unit outward normal vectors on $\partial \Omega _{S}$ and $\partial \Omega _{D}$, respectively, and $\mathbf{\tau }_{j}\ (j=1,\cdots ,d-1)$ as the unit tangential vectors on the interface $\Gamma $. We note that  
 $\mathbf{n}_{S}=-\mathbf{n}_{D}~\mathrm{on}~\Gamma $, and the notations  $\Gamma_{S}=\partial \Omega _{S}\setminus \Gamma ,\Gamma _{D}=\partial \Omega_{D}\setminus \Gamma $ are also used, see Fig. 2.1 for a sketch of the problem domain setting.

\begin{figure}[htbp]
\centering
\includegraphics[width=90mm,height=46mm]{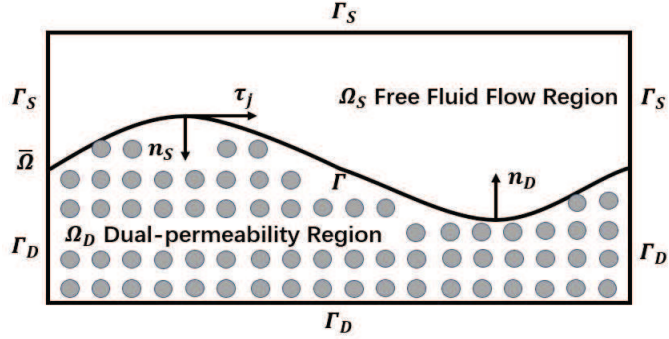}
\caption{The global domain $\overline{\Omega}$ consisting of the fluid region
$\Omega_S$ and the dual-permeability region $\Omega_D$ separated by the
interface $\Gamma$.}
\end{figure}

In the fluid region $\Omega_{S} $, the fluid velocity $\mathbf{u}_{S}$ and kinematic pressure $p_{S}$
are assumed to satisfy the Stokes equations:
\begin{eqnarray}
 -\nabla \cdot \mathbb{T}(\mathbf{u}_{S},p_S)&=&\mathbf{f}_{S}\ \ \ \ \mathrm{in}~\Omega _{S},
\label{Stokes1} \\
\nabla \cdot \mathbf{u}_{S} &=&0\ \ \ \ \hspace{1.5mm} \mathrm{in}~\Omega _{S},\label{Stokes2}
\end{eqnarray}
where $\mathbb{T}(\mathbf{u}_{S},p_S)=-p_s\mathbb{I}+2\nu\mathbb{D}(\mathbf{u}_{S})$ indicates the stress tensor, among that $\mathbb{I}$ is an identity matrix and $\mathbb{D}( \mathbf{u}_{S})=\frac{1}{2}(\nabla \mathbf{u}_{S}+(\nabla \mathbf{u}_{S})^{T})$ means the deformation tensor, $\nu $ represents the kinematic viscosity of the fluid flow. Besides $\mathbf{f}_S$ expresses the external body force.

In the dual-permeability region $\Omega_{D} $, we  utilize the traditional dual-permeability model, which is composed of microfracture flow and matrix flow equations \cite{HouJY16, Mahbub2019}. The microfracture flow is governed by the following mixed equations with the microfracture flow velocity $\mathbf{u}_{D}$ and the microfracture flow pressure $\varphi_{D}$:
\begin{eqnarray}
\mathbf{u}_{D}&=&-\frac{k_D}{\mu} \nabla \varphi_{D} \hspace{8mm}
\mathrm{in}~\Omega_{D} ,  \label{Darcy1} \\
 \nabla \cdot \mathbf{u}_{D} + \frac{\sigma k_M}{\mu} (\varphi_D-\varphi_M) &=& {f}_{D}
 \hspace{20mm} \mathrm{in}~\Omega _{D}. \label{Darcy2}
\end{eqnarray}
In order to describe the matrix flow velocity $\mathbf{u}_{M}$ and the matrix flow pressure $\varphi_{M}$, we can adopt the following system of equations:
\begin{eqnarray}
\mathbf{u}_{M}&=&-\frac{k_M}{\mu} \nabla \varphi_{M} \hspace{6.4mm}
\mathrm{in}~\Omega_{D} ,  \label{Darcy3} \\
 \nabla \cdot \mathbf{u}_{M} + \frac{\sigma k_M}{\mu} (\varphi_M-\varphi_D) &=& {f}_{M}
 \hspace{20.6mm} \mathrm{in}~\Omega _{D}. \label{Darcy4}
\end{eqnarray}
Here, $\sigma$ indicates the shape factor which characterizes the morphology and dimension between the microfractures and matrix. $\mu$ is the dynamic viscosity of the dual-permeability flow. $k_D$ and $k_M$ are the intrinsic permeabilities in microfractures and matrix respectively. In addition, ${f}_{D}$ and $f_M $ denote the sink/source terms (here $f_M $ is usually considered to be equal to zero). Note that $f_D$ needs to satisfy the solvability condition $\int_{\Omega_D} f_D=0$. As for the term $\frac{\sigma k_M}{\mu} (\varphi_D-\varphi_M)$, it means the mass exchange between microfractures and matrix.

We assume the impermeable boundary conditions in the sense that the fluid velocity $\mathbf{u}_{S} $, the microfracture flow velocity $\mathbf{u}_{D}$ and  the matrix flow velocity $\mathbf{u}_{M}$ satisfy the homogeneous Dirichlet boundary conditions, i.e., no slip condition $\mathbf{u}_{S}=0~~\mathrm{on}~\Gamma _{S} $ in the Stokes region, and impermeable conditions
$\mathbf{u}_{D}\cdot\mathbf{n}_D=0, ~ \mathbf{u}_{M}\cdot\mathbf{n}_D=0~~\mathrm{on}~\Gamma _{D} $ in the dual-permeability region.
Such boundary conditions can be easily extended.

As to the interface coupling conditions on $\Gamma$, similarly to the requirement in traditional Stokes-Darcy model, the conservation of mass, the balance of forces and the tangential conditions for the velocity need to be enforced. In this paper, we adopt the important Beavers-Joseph (BJ) interface conditions as the tangential conditions on the interface $\Gamma $, see \cite{Beavers,Jones}. To this end, the interface coupling conditions between the flow in the microfractures and the flow in the conduits/macrofractures are assumed as follows:
\begin{eqnarray}
\mathbf{u}_{S}\cdot \mathbf{n}_{S}+\mathbf{u}_{D}\cdot
\mathbf{n}_{D}
&=&0\ \ \ \ \ \ \ \ \ \ \ \ \ \ \hspace{7.5mm} \mathrm{on}~\Gamma ,  \label{interface1} \\
-\mathbf{n}_{S}\cdot (\mathbb{T}(\mathbf{u}_{S},p_S)\cdot \mathbf{n}_{S})
&=&  \frac{1}{\rho}\varphi_{D} \hspace{18.3mm} \mathrm{on}~\Gamma, \label{interface2}\\
-\tau_{j}\cdot(\mathbb{T}(\mathbf{u}_{S},p_S)\cdot \mathbf{n}_{S})
&=&\frac{\nu\alpha\sqrt{d}}{\sqrt{\mathrm{trace}(\mathbf{\prod})}} \tau _{j} \cdot (\mathbf{u}_{S}-\mathbf{u}_{D})
\hspace{3mm}  1\leq j \leq d-1\hspace{3mm} \mathrm{on}~\Gamma,\label{interface3BJ}
\end{eqnarray}
where $\alpha$ represents an experimentally determined positive parameter depending on the porous medium properties, $\rho$ is the fluid density and $\mathbf{\prod} = k_D \mathbb{I}$.
Moreover, Hou et al. \cite{HouJY16} proposed a no-exchange condition between the matrix and the conduits/macrofractures, which means no flux could pass across the interface. The following condition is used as the fourth interface condition:
\begin{eqnarray}
\mathbf{u}_{M}\cdot \mathbf{n}_{D} &=&0\ \ \ \ \ \ \ \ \ \ \ \ \ \ \hspace{7.5mm} \mathrm{on}~\Gamma.  \label{interface4}
\end{eqnarray}

In order to solve the coupled Stokes-dual-permeability model by DDM, a natural idea is to consider Robin-type boundary conditions for both the Stokes equations and the dual-permeability system respectively. Firstly for the Stokes equations, two Robin-type functions $g_S,\ g_{S,\tau}$  on $\Gamma$ are constructed: for a given constant $\delta_S>0$, 
\begin{eqnarray}
&&g_S=-\mathbf{n}_{S}\cdot (\mathbb{T}(\mathbf{u}_S,p_S)\cdot \mathbf{n}_{S})
-\delta_S{\mathbf{u}}_{S}\cdot \mathbf{n}_{S},\label{Robintype1} \\
&&g_{S,\tau}=\sum_{j=1}^{d-1}g_{S,\tau_j}=\sum_{j=1}^{d-1}\Big(-\tau_{j}\cdot \mathbb{T}(\mathbf{u}_S,p_S)\cdot \mathbf{n}_{S}
-\frac{\nu\alpha\sqrt{d}}{\sqrt{\mathrm{trace}(\mathbf{\prod})}}\hspace{0.5mm}\mathbf{u}_S \cdot \tau _{j}\Big).\label{Robintype2}
\end{eqnarray}
Similarly, we propose the following Robin-type condition for the dual-permeability system in the porous media: for a given constant $\delta_D>0$, define the Robin-type function $g_D$ on $\Gamma$ as
\begin{eqnarray}
g_D&=&\frac{1}{\rho}\varphi_D-\delta_D\mathbf{u}_D \cdot\mathbf{n}_D.\label{Robintype3}
\end{eqnarray}

In the following Lemma, we present the equivalence of the original interface conditions
(\ref{interface1})-(\ref{interface3BJ}) and the above Robin-type
conditions (\ref{Robintype1})-(\ref{Robintype3}).

\begin{lemma}
The interface coupling conditions
(\ref{interface1})-(\ref{interface3BJ}) is equivalent to the three Robin-type
conditions (\ref{Robintype1})-(\ref{Robintype3}) if and
only if $g_S$, $g_{S,\tau}$, and $g_D$ satisfy the following three compatibility conditions
on the interface $\Gamma$:
\begin{eqnarray}
&&g_D =g_S+(\delta_S+\delta_D)
 {\mathbf{u}}_{S}\cdot
\mathbf{n}_{S}, \label{comp1}\\
&&g_S =g_D+(\delta_S+\delta_D)
{\mathbf{u}}_{D}\cdot \mathbf{n}_{D},\label{comp2}\\
&&g_{S,\tau}=\sum_{j=1}^{d-1}g_{S,\tau_j}=\sum_{j=1}^{d-1}\Big(-\frac{\nu\alpha\sqrt{d}}{\sqrt{\mathrm{trace}(\mathbf{\prod})}}
\mathbf{u}_D \cdot \tau_{j}\Big).\label{comp3}
\end{eqnarray}
\end{lemma}
\begin{proof}
Adding the compatibility conditions (\ref{comp1})-(\ref{comp2}) together, the first original interface condition (\ref{interface1}) can be demonstrated directly. Substituting the third compatibility condition (\ref{comp3}) into (\ref{Robintype2}), we have
\begin{eqnarray*}
-\frac{\nu\alpha\sqrt{d}}{\sqrt{\mathrm{trace}(\mathbf{\prod})}}\sum_{j=1}^{d-1}\mathbf{u}_D\cdot \tau _{j}=-\sum_{j=1}^{d-1}(\tau_{j}\cdot \mathbb{T}(\mathbf{u}_S,p_S)\cdot \mathbf{n}_{S})
-\frac{\nu\alpha\sqrt{d}}{\sqrt{\mathrm{trace}(\mathbf{\prod})}}\sum_{j=1}^{d-1}\mathbf{u}_S \cdot \tau _{j}.
\end{eqnarray*}
Hence the original interface condition (\ref{interface3BJ}) has been proved. By combining (\ref{comp2}) with (\ref{Robintype3}), it follows that
\begin{eqnarray}\label{gS}
g_S=\frac{1}{\rho}\varphi_D-\delta_D\mathbf{u}_D \cdot \mathbf{n}_D
+(\delta_S+\delta_D){\mathbf{u}}_{D}\cdot \mathbf{n}_{D}.
\end{eqnarray}
Substituting (\ref{gS}) into the Robin-type condition (\ref{Robintype1}), and utilizing the equation (\ref{interface1}), the second original interface coupling condition (\ref{interface2}) can be directly obtained. So we drive the original interface coupling conditions (\ref{interface1})-(\ref{interface3BJ}) by the Robin-type conditions (\ref{Robintype1})-(\ref{Robintype3}) and the compatibility conditions (\ref{comp1})-(\ref{comp3}) on the interface $\Gamma$. Similarly, we can also demonstrate the Robin-type conditions (\ref{Robintype1})-(\ref{Robintype3}) by (\ref{interface1})-(\ref{interface3BJ}) and (\ref{comp1})-(\ref{comp3}), which means the equivalence of two types of interface conditions under the compatibility conditions (\ref{comp1})-(\ref{comp3}).

Conversely, Suppose that (\ref{interface1})-(\ref{interface3BJ}) and (\ref{Robintype1})-(\ref{Robintype3}) are equivalent, we can verify the necessary compatibility
conditions  (\ref{comp1})-(\ref{comp3}) for $g_{S}, g_{S,\tau} $ and $\ g_{D}$. Actually, combining (\ref{interface3BJ}) with (\ref{Robintype2}), we can immediately arrive at the third compatibility condition (\ref{comp3}). By substituting (\ref{interface2}) into (\ref{Robintype1}) and then applying (\ref{interface1}) and (\ref{Robintype3}), we can deduce:
\begin{eqnarray*}
	g_{S}&=&\frac{1}{\rho}\varphi_{D}+\delta_S\mathbf{u}_{D}\cdot\mathbf{n}_D\\
	&=&\frac{1}{\rho}\varphi_{D}-\delta_D\mathbf{u}_{D}\cdot\mathbf{n}_D+\delta_D\mathbf{u}_{D}\cdot\mathbf{n}_D+\delta_S\mathbf{u}_{D}\cdot\mathbf{n}_D\\
	&=&g_{D}+(\delta_S+\delta_D)\mathbf{u}_{D}\cdot\mathbf{n}_D,
\end{eqnarray*}
which means that the second compatibility condition (\ref{comp2}) holds. Condition (\ref{comp1}) can be verified by very similar argument.
\end{proof}

\section{ The Modified Weak Formulation and its Equivalence }
Firstly, we introduce some Sobolev spaces and norms. For the fluid domain $\Omega_S $ and the porous media domain $\Omega_D $, we denote the inner products by $(\cdot, \cdot )_{S}$ and $(\cdot, \cdot )_{D}$ respectively, and the corresponding $L^2$-norms by $||\cdot ||_S$ and $||\cdot ||_D$. Meanwhile, $\langle \cdot, \cdot \rangle$ is defined as the $L^2$ inner product on the interface $\Gamma$, and the $L^{2}(\Gamma )$-norm is denoted by $||\cdot ||_{\Gamma }$. By setting the space
\begin{eqnarray*}
H_{\mathrm{div}}=H(\mathrm{div};\Omega _{D}):=\{\mathbf{v}_{D}\in L^{2}(\Omega
_{D})^{d}:\nabla \cdot \mathbf{v}_{D}\in L^{2}(\Omega _{D})\},
\end{eqnarray*}
the Sobolev spaces for the Stokes-dual-permeability model are defined as follows:
\begin{eqnarray*}
\mathbf{X}^0_{S}&:=&\{\mathbf{v}_{S}\in H^{1}(\Omega _{S})^{d}:\mathbf{v}_{S}=0~%
\mathrm{on}~\Gamma _{S}\},\ \ \hspace{19mm}
Q^0_{S}:=L^{2}(\Omega _{S}), \\
\mathbf{X}^0_{D}&:=&\{\mathbf{v}_{D}\in H(\mathrm{div};\Omega
_{D}):\mathbf{v}_{D}\cdot
\mathbf{n}_{D}=0~%
\mathrm{on}~\Gamma _{D}\}, \ \ \hspace{6mm} Q^0_{D}:=L^{2}(\Omega _{D}),\\
\mathbf{X}^0_{M}&:=&\{\mathbf{v}_{M}\in H(\mathrm{div};\Omega
_{D}):\mathbf{v}_{M}\cdot
\mathbf{n}_{D}=0~%
\mathrm{on}~\partial\Omega_D\}, \hspace{5mm} Q^0_{M}:=L^{2}(\Omega _{D}).
\end{eqnarray*}

Then we can construct two product spaces as 
\begin{eqnarray*}
\mathbf{X}^0&:=&\Big\{\mathbf{v}=(\mathbf{v}_S,\mathbf{v}_D,\mathbf{v}_M)
\in\mathbf{X}^0_S\times\mathbf{X}^0_D\times\mathbf{X}^0_M\Big\},\\
Q^0&:=&\Big\{q=(q_S,\psi_D,\psi_M)\in Q^0_S\times Q^0_D\times Q^0_M:
\int_{\Omega_S} q_S + \int_{\Omega_D} \psi_D =0\Big\}.
\end{eqnarray*}
Based on space settings and the Robin-type boundary conditions defined above, we firstly introduce the classical weak formulation for the coupled model in the the global domain $\Omega$,  as follows.

\textbf{Original weak formulation.} The weak formulation for the coupled Stokes-dual-permeability
problem with the Robin-type boundary conditions (\ref{Robintype1})-(\ref{Robintype3}) is formulated as follows: for given $g_S, \ g_{S,\tau}, \ g_D \in L^2(\Gamma)$, find $\mathbf{u}^0:=({\mathbf{u}}^0_{S},{\mathbf{u}}^0_{D},{\mathbf{u}}^0_{M})\in \mathbf{X}^0$ and $p^0:=(p^0_S,\varphi^0_D,\varphi^0_M) \in Q^0$, such that
\begin{eqnarray}\label{weak-o1}
a(\mathbf{u}^0,\mathbf{v})-b(\mathbf{v},{p}^0)
+b_{DM}(\mathbf{v},p^0)
+c_{\Gamma}(\mathbf{u}^0, \mathbf{v})
&=&L(\mathbf{v})
\hspace{15mm} \forall ~\mathbf{v}:=(\mathbf{v}_S,\mathbf{v}_D,\mathbf{v}_M)\in \mathbf{X}^0, \\
\label{weak-o2}b(\mathbf{u}^0,q)+a_{DM}(p^0,q)&=&  \frac{1}{\rho}(f_{D},\psi_D)_D \hspace{4mm} \forall ~q:=(q_S,\psi_D,\psi_M)\in Q^0,
\end{eqnarray}
with the compatibility conditions (\ref{comp1})-(\ref{comp3}) on $\Gamma$. Here, the following bilinear forms are used:
\begin{eqnarray*}
a(\mathbf{u},\mathbf{v})&:=& a_{S}(\mathbf{u}_{S},\mathbf{v}_{S})+a_{D}(\mathbf{u}_{D},\mathbf{v}_{D})+a_{M}(\mathbf{u}_{M},\mathbf{v}_{M}),\\
b(\mathbf{v},p)&:=&b_{S}(\mathbf{v}_{S},p_S)+b_{D}(\mathbf{v}_{D},\varphi_D )+b_{M}(\mathbf{v}_{M},\varphi_M ),\\
b_{DM}(\mathbf{v},p)&:=& \frac{\sigma k_M}{\rho k_D}(\varphi_D-\varphi_M, \nabla\cdot\mathbf{v}_D)_D
+\frac{\sigma }{\rho }(\varphi_M-\varphi_D, \nabla\cdot\mathbf{v}_M)_D,\\
c_{\Gamma}(\mathbf{u}, \mathbf{v})&:=&
\delta_S\langle{\mathbf{u}}_{S}\cdot \mathbf{n}_S,\mathbf{v}_{S}\cdot \mathbf{n}_S\rangle
+\delta_D\langle{\mathbf{u}}_{D}\cdot \mathbf{n}_D,\mathbf{v}_{D}\cdot \mathbf{n}_D\rangle,\\
L(\mathbf{v})&:=&-\langle{g}_{S},\mathbf{v}_{S}\cdot\mathbf{n}_S\rangle
-\sum_{j=1}^{d-1}\langle g_{S,\tau_j},\mathbf{v}_D\cdot\tau_j \rangle
+(\mathbf{f}_{S},\mathbf{v}_{S})_S-\langle{g}_{D},\mathbf{v}_{D}\cdot \mathbf{n}_D\rangle
+ \frac{\mu}{\rho k_D}(f_{D},\nabla\cdot \mathbf{v}_{D})_D,\\
a_{DM}(p,q)&:=& \frac{\sigma k_M}{\rho \mu}(\varphi_D-\varphi_M, \psi_D)_D
+\frac{\sigma k_M}{\rho \mu}(\varphi_M-\varphi_D, \psi_M)_D,
\end{eqnarray*}
with
\begin{eqnarray*}
a_{S}(\mathbf{u}_{S},\mathbf{v}_{S})&:=& 2\nu (\mathbb{D}(\mathbf{u}_{S}),\mathbb{D}(\mathbf{v}_{S}))_S
+\sum_{j=1}^{d-1}\frac{\nu\alpha\sqrt{d}}{\sqrt{\mathrm{trace}(\prod)}}
\langle\mathbf{u}_{S}\cdot \tau _{j},\mathbf{v}_{S}\cdot \tau _{j}\rangle,\quad b_{S}(\mathbf{v}_{S},p_S):=(p_S,\nabla \cdot \mathbf{v}_{S})_S,\\
a_{D}(\mathbf{u}_{D},\mathbf{v}_{D})&:=& \frac{\mu}{\rho k_D}(\mathbf{u}_{D},\mathbf{v}_{D})_D
+\frac{\mu}{\rho k_D}(\nabla\cdot \mathbf{u}_{D},\nabla\cdot \mathbf{v}_{D})_D,\quad b_{D}(\mathbf{v}_{D},\varphi_D ):=  \frac{1}{\rho}(\varphi_D ,\nabla \cdot\mathbf{v}_{D})_D,\\
a_{M}(\mathbf{u}_{M},\mathbf{v}_{M})&:=& \frac{\mu}{\rho k_M}(\mathbf{u}_{M},\mathbf{v}_{M})_D
+\frac{\mu }{\rho k_M}(\nabla\cdot \mathbf{u}_{M},\nabla\cdot \mathbf{v}_{M})_D,\quad
b_{M}(\mathbf{v}_{M},\varphi_M ):=  \frac{1}{\rho}(\varphi_M ,\nabla \cdot\mathbf{v}_{M})_D.
\end{eqnarray*}

In order to decouple the original Stokes-dual-permeability problem into the independent Stokes subproblem and dual-permeability subproblems, a modified meaningful weak formulation will be constructed. Moreover, the convergence analysis of our algorithm in the following sections is based on this modified weak formulation, which is inspired by Galvis et al. \cite{Twoweakform} and Sun et al. \cite{Sun2021}. We define some modified Sobolev spaces as follows,
\begin{eqnarray*}
\mathbf{X}_{S}&:=&\Big\{\mathbf{v}_{S}\in \mathbf{X}_S^0: \int_{\Gamma}\mathbf{v}_{S} \cdot \mathbf{n}_S=0 ~\mathrm{on}~\Gamma\Big\},\hspace{12mm}
Q_{S}:=L_0^{2}(\Omega _{S}), \\
\mathbf{X}_{D}&:=&\Big\{\mathbf{v}_{D}\in \mathbf{X}_D^0: \int_{\Gamma}\mathbf{v}_{D} \cdot \mathbf{n}_D=0 ~\mathrm{on}~\Gamma \Big\}, \hspace{9.5mm} Q_{D}:=L_0^{2}(\Omega _{D}),\\
\mathbf{X}_{M}&:=& \mathbf{X}_M^0, \hspace{56.5mm} Q_{M}:=L_0^{2}(\Omega _{D}),
\end{eqnarray*}%
with the following norms,
\begin{eqnarray*}
||\mathbf{v}_{S}||_{1}&=&\sqrt{||\mathbf{v}_{S}||_S^{2}+||\nabla \mathbf{v}_{S}||_S^{2}} \ \hspace{20.7mm} \ \ \forall\  \mathbf{v}_{S}\in \mathbf{X}_{S}, \\
||\mathbf{v}_{D}|| _{\mathrm{div}}&=&\sqrt{||\mathbf{v}_{D}||_D^{2}+
||\nabla \cdot \mathbf{v}_{D}||_D^{2}}\ \ \ \hspace{1mm}
\ \ \ \ \ \ \ \ \ \ \ \ \ \forall\  \mathbf{v}_{D}\in
\mathbf{X}_{D},\\
||\mathbf{v}_{M}|| _{\mathrm{div}}&=&\sqrt{||\mathbf{v}_{M}||_D^{2}+
||\nabla \cdot \mathbf{v}_{M}||_D^{2}}\ \ \ \hspace{1mm}
\ \ \ \ \ \ \ \ \ \ \ \ \forall\  \mathbf{v}_{M}\in
\mathbf{X}_{M}.
\end{eqnarray*}
Furthermore, we can define the product spaces $\mathbf{X}:=\mathbf{X}_S \times \mathbf{X}_D \times \mathbf{X}_M$ and $Q:=Q_S \times Q_D \times Q_M$ in the global domain $\Omega$.

\textbf{Modified weak formulation.} The modified meaningful weak formulation for the coupled Stokes-dual-permeability system on the modified Sobolev spaces with the Robin-type boundary conditions (\ref{Robintype1})-(\ref{Robintype3}) is introduced as: for given $g_S, \ g_{S,\tau}, \ g_D \in L^2(\Gamma)$, find $\mathbf{u}:=({\mathbf{u}}_{S},{\mathbf{u}}_{D},{\mathbf{u}}_{M})\in \mathbf{X}$ and $p:=(p_S,\varphi_D,\varphi_M) \in Q$, such that
\begin{eqnarray}\label{weak-m1}
a(\mathbf{u},\mathbf{v})-b(\mathbf{v},{p})
+b_{DM}(\mathbf{v},p)
+c_{\Gamma}(\mathbf{u}, \mathbf{v})
&=&L(\mathbf{v})
\hspace{15mm} \forall ~\mathbf{v}:=(\mathbf{v}_S,\mathbf{v}_D,\mathbf{v}_M)\in \mathbf{X}, \\
\label{weak-m2}b(\mathbf{u},q)+a_{DM}(p,q)&=&  \frac{1}{\rho}(f_{D},\psi_D)_D \hspace{4mm} \forall ~q:=(q_S,\psi_D,\psi_M)\in Q,
\end{eqnarray}
with the compatibility conditions (\ref{comp1})-(\ref{comp3}) on the interface $\Gamma$.

\begin{theorem}\label{WFequal}
The original weak formulation (\ref{weak-o1})-(\ref{weak-o2}) is equivalent to the modified weak formulation (\ref{weak-m1})-(\ref{weak-m2}) 
in the following sense:
\begin{enumerate}
\item Two weak formulations have the same velocity solutions.
\item If we get pressure solutions of one weak formulation, the solution of another weak formulation can be transformed with the mean value of pressure in the global domain $\Omega$.
\end{enumerate}
\end{theorem}

\begin{proof}
($\Rightarrow$) Suppose that the solutions of the original weak formulation (\ref{weak-o1})-(\ref{weak-o2}) are $\mathbf{u}^0=(\mathbf{u}^0_S,\mathbf{u}^0_D,\mathbf{u}^0_M)\in\mathbf{X}^0$ and $p^0=(p_S^0,\varphi^0_D, \varphi^0_M)\in Q^0$. By the compatibility conditions (\ref{comp1})-(\ref{comp2}), we can yield
$\mathbf{u}^0_S\cdot \mathbf{n}_S + \mathbf{u}^0_D\cdot \mathbf{n}_D=0$ on $\Gamma$, so that
\begin{eqnarray*}
\int_{\Gamma}\mathbf{u}^0_S\cdot \mathbf{n}_S + \int_{\Gamma} \mathbf{u}^0_D\cdot \mathbf{n}_D=0.
\end{eqnarray*}
Introducing the regional average operator $p^c=(\frac{1}{|\Omega_S|},-\frac{\rho}{|\Omega_D|}, -\frac{\rho}{|\Omega_D|})\in Q^0$, where $|\Omega_{S(D)}|=\int_{\Omega_{S(D)}}1$, we can get $a_{DM}(p^0,p^c)= 0$. Then applying the solvability condition $\int_{\Omega_D}f_D=0$, we can have
\begin{eqnarray*}
0=b(\mathbf{u}^0,p^c)&=&\frac{1}{|\Omega_S|}\int_{\Omega_S}\nabla \cdot\mathbf{u}^0_S
- \frac{1}{|\Omega_D|}\int_{\Omega_D} \nabla \cdot \mathbf{u}^0_D
- \frac{1}{|\Omega_D|}\int_{\Omega_D} \nabla \cdot \mathbf{u}^0_M\\
&=&\frac{1}{|\Omega_S|}\int_{\Gamma}\mathbf{u}^0_S\cdot \mathbf{n}_S
- \frac{1}{|\Omega_D|} \int_{\Gamma} \mathbf{u}^0_D\cdot \mathbf{n}_D.
\end{eqnarray*}
We can directly obtain $\int_{\Gamma}\mathbf{u}^0_S \cdot \mathbf{n}_S =\int_{\Gamma}\mathbf{u}^0_D \cdot \mathbf{n}_D =0 $, which means $\mathbf{u}^0\in \mathbf{X}$. For the pressure, we construct the solution $p\in Q$ by removing the mean value of $p^0$ as follows,
\begin{eqnarray*}
p:=\Big(p^0_S-\frac{1}{|\Omega_S|}\int_{\Omega_S}p^0_S,\hspace{1mm}
\varphi^0_D-\frac{\rho}{|\Omega_D|}\int_{\Omega_D}\varphi^0_D,\hspace{1mm}
\varphi^0_M-\frac{\rho}{|\Omega_D|}\int_{\Omega_D}\varphi^0_M\Big)\in Q.
\end{eqnarray*}
So it is clearly to arrive at $b(\mathbf{v},p)=b(\mathbf{v},p^0)$ for all $\mathbf{v}\in \mathbf{X}$. We also can get
\begin{eqnarray*}
b_{DM}(\mathbf{v},p) &=& \frac{\sigma k_M}{\rho k_D}
\Big(\varphi^0_D-\frac{1}{|\Omega_D|}\int_{\Omega_D}\varphi^0_D
-\varphi^0_M+\frac{1}{|\Omega_D|}\int_{\Omega_D}\varphi^0_M,\hspace{1mm} \nabla\cdot\mathbf{v}_D\Big)_D\\
&&+\frac{\sigma }{\rho } \Big(\varphi^0_M-\frac{1}{|\Omega_D|}\int_{\Omega_D}\varphi^0_M
-\varphi^0_D+\frac{1}{|\Omega_D|}\int_{\Omega_D}\varphi^0_D, \hspace{1mm} \nabla\cdot\mathbf{v}_M\Big)_D\\
&=& \frac{\sigma k_M}{\rho k_D}(\varphi^0_D-\varphi^0_M,\nabla\cdot\mathbf{v}_D)_D
+\frac{\sigma }{\rho } (\varphi^0_M-\varphi^0_D,\nabla\cdot\mathbf{v}_M)_D = b_{DM}(\mathbf{v},p^0).
\end{eqnarray*}
Moreover, for all $q \in Q$, we can yield
\begin{eqnarray*}
a_{DM}(p,q) = \frac{\sigma k_M}{\rho \mu}(\varphi^0_D-\varphi^0_M, \psi_D)_D
+\frac{\sigma k_M}{\rho \mu}(\varphi^0_M-\varphi^0_D, \psi_M)_D = a_{DM}(p^0,q).
\end{eqnarray*}
Hence, we can summarize that ($\mathbf{u}^0, p$) is the solution of the modified weak formulation (\ref{weak-m1})-(\ref{weak-m2}).

($\Leftarrow$) Let $\mathbf{u}=(\mathbf{u}_S,\mathbf{u}_D, \mathbf{u}_M)\in\mathbf{X}$ and $p=(p_S,\varphi_D,\varphi_M)\in Q$ be the solutions of the original weak formulation (\ref{weak-o1})-(\ref{weak-o2}). Set $\mathbf{w}=(\mathbf{w}_S,\mathbf{w}_D,\mathbf{w}_M)\in \mathbf{X}^0$ such that
\begin{eqnarray*}
\mathbf{w}_S\cdot \mathbf{n}_S=\frac{|\Omega_S|}{2|\Gamma|} \hspace{3mm} \mathrm{on} ~\Gamma,
\hspace{6mm} \mathbf{w}_D \cdot \mathbf{n}_D= -\frac{|\Omega_D|}{2|\Gamma|}
 \hspace{3mm} \mathrm{on}~ \Gamma, \hspace{6mm} \mathbf{w}_M \cdot \mathbf{n}_D= 0
 \hspace{3mm} \mathrm{on}~ \Gamma.
\end{eqnarray*}
Recalling that the regional average operator $p^c=(\frac{1}{|\Omega_S|},-\frac{\rho}{|\Omega_D|},-\frac{\rho}{|\Omega_D|}) \in Q^0$, we get
\begin{eqnarray*}
b(\mathbf{w},p^c)&=&\frac{1}{|\Omega_S|}\int_{\Omega_S}\nabla \cdot\mathbf{w}_S
- \frac{1}{|\Omega_D|}\int_{\Omega_D} \nabla \cdot \mathbf{w}_D
- \frac{1}{|\Omega_D|}\int_{\Omega_D} \nabla \cdot \mathbf{w}_M\\
&=& \frac{1}{|\Omega_S|} \int_{\Gamma}\mathbf{w}_S\cdot \mathbf{n}_S
- \frac{1}{|\Omega_D|} \int_{\Gamma} \mathbf{w}_D\cdot \mathbf{n}_D=1.
\end{eqnarray*}
Then we can define $p^0=(p+\beta p^c) \in Q^0$, where
\begin{eqnarray*}
\beta:=a(\mathbf{u},\mathbf{w})-b(\mathbf{w},p)+b_{DM}(\mathbf{w},p)+c_{\Gamma}(\mathbf{u},\mathbf{w})
-L(\mathbf{w}).
\end{eqnarray*}
From now on, we need to demonstrate that $\mathbf{u}=(\mathbf{u}_S,\mathbf{u}_D,\mathbf{u}_M)$ and $p^0=(p_S^0,\varphi^0_D,\varphi^0_M)$ are the solutions of the original weak formulation (\ref{weak-o1})-(\ref{weak-o2}). For any $(\mathbf{v}^0,q^0)\in(\mathbf{X^0},Q^0)$, there exists a constant $\chi$ such that $\mathbf{v}:=(\mathbf{v}^0+\chi\mathbf{w}) \in \mathbf{X}$. Setting $q=(q_S,\psi_D,\psi_M):=\Big(q^0_S-\frac{1}{|\Omega_S|}\int_{\Omega_S}q^0_S, \hspace{1mm}
\psi^0_D-\frac{\rho}{|\Omega_D|}\int_{\Omega_D}\psi^0_D, \hspace{1mm} \psi^0_M-\frac{\rho}{|\Omega_D|}\int_{\Omega_D}\psi^0_M\Big) \in Q$, we can arrive at
\begin{eqnarray*}
a(\mathbf{u},\mathbf{v}^0)&-&b(\mathbf{v}^0,p^0)+b_{DM}(\mathbf{v}^0,p^0)+c_{\Gamma}(\mathbf{u},\mathbf{v}^0)\\
&=&\Big\{a(\mathbf{u},\mathbf{v})-b(\mathbf{v},p)+b_{DM}(\mathbf{v},p)+c_{\Gamma}(\mathbf{u},\mathbf{v}) \Big\}-\beta b(\mathbf{v},p^c)-\beta b_{DM}(\mathbf{v},p^c)\\
&&-\chi\Big\{a(\mathbf{u},\mathbf{w})-b(\mathbf{w},p)+b_{DM}(\mathbf{w},p)
+c_{\Gamma}(\mathbf{u},\mathbf{w})-\beta b(\mathbf{w},p^c)-\beta b_{DM}(\mathbf{w},p^c) \Big\}\\
&=&L(\mathbf{v})-\chi L(\mathbf{w})=L(\mathbf{v}^0),\\
b(\mathbf{u},q^0) + a_{DM}(p^0,q^0) &=& b(\mathbf{u},q) + a_{DM}(p,q) + \beta a_{DM}(p^c,q^0)\\
&=&  b(\mathbf{u},q) + a_{DM}(p,q) = \frac{1}{\rho}(f_D,\psi_D)_D = \frac{1}{\rho}(f_D,\psi_D^0)_D.
\end{eqnarray*}
Overall, we have proved the equivalence of two weak formulations.
\end{proof}

\textbf{Decoupled modified weak formulation.} The weak formulation of the decoupled Stokes and dual-permeability
model with Robin-type boundary conditions (\ref{Robintype1})-(\ref{Robintype3}) is given as follows: for two given functions $g_S,\ g_D\in L^2(\Gamma)$, find $({\mathbf{u}}_{S},{p_S}; {\mathbf{u}}_{D},{\varphi_D}; {\mathbf{u}}_{M},{\varphi_M}) \in  (\mathbf{X}_{S}, Q_{S}; \mathbf{X}_{D},Q_{D}; \mathbf{X}_{M},Q_{M})$,
such that, for any $(\mathbf{v}_S,q_S;\mathbf{v}_D,\psi_D; \mathbf{v}_M,\psi_M)\in (\mathbf{X}_{S},Q_S;\mathbf{X}_D,Q_D;\mathbf{X}_M,Q_M)$, we have
\begin{eqnarray}
a_{S}(\mathbf{u}_{S},\mathbf{v}_{S})-b_S(\mathbf{v}_S, p_S)
+\delta_S\langle\mathbf{u}_S\cdot\mathbf{n}_S,&\mathbf{v}_{S}&\cdot\mathbf{n}_S\rangle \nonumber\\
&=&(\mathbf{f}_S,\mathbf{v}_S)_{S}-\langle g_S,\mathbf{v}_S\cdot\mathbf{n}_S \rangle
-\sum_{j=1}^{d-1}\langle g_{S,\tau_j},\mathbf{v}_S\cdot \tau_{j} \rangle,\label{StokesR1} \\
b_S(\mathbf{u}_S, q_S)&=&0,\label{StokesR2}\\
a_{D}(\mathbf{u}_{D},\mathbf{v}_{D})-b_D(\mathbf{v}_D,\varphi_D)
+\frac{\sigma k_M}{\rho k_D}(\varphi_D&-&\varphi_M, \nabla \cdot \mathbf{v}_D)_D
+\delta_D\langle\mathbf{u}_D\cdot\mathbf{n}_D,\mathbf{v}_D\cdot\mathbf{n}_D\rangle\nonumber\\
&=&\frac{\mu}{\rho k_D}(f_D,\mathrm{div} \hspace{0.5mm} \mathbf{v}_D)_{D}
-\langle g_D,\mathbf{v}_D\cdot\mathbf{n}_D\rangle,\label{DarcyR1} \\
b_D(\mathbf{u}_D,\psi_D)
+\frac{\sigma k_M}{\rho \mu}(\varphi_D-\varphi_M, \psi_D)_D
&=&\frac{1}{\rho}(f_D,\psi_D)_{D},\label{DarcyR2}\\
a_{M}(\mathbf{u}_{M},\mathbf{v}_{M})-b_M(\mathbf{v}_M,\varphi_M)
&+&\frac{\sigma }{\rho }(\varphi_M-\varphi_D, \nabla \cdot \mathbf{v}_M)_D = 0,\label{DarcyR3} \\
b_M(\mathbf{u}_M,\psi_M)
+\frac{\sigma k_M}{\rho \mu}(\varphi_M-\varphi_D, \psi_M)_D
&=& 0.\label{DarcyR4}
\end{eqnarray}
with the compatibility conditions (\ref{comp1})-(\ref{comp3}) on the interface $\Gamma$.
It is worth to mention that the well-posedness for the above weak formulation of the Stokes-dual-permeability decoupled model can be certified clearly.

\section{ Parallel Robin-type Domain Decomposition Method}
Now, we propose the following parallel Robin-type DDM for solving the fully-mixed coupled Stokes-dual-permeability  problem with the BJ interface conditions.

\textbf{{{Parallel DDM Algorithm}} (PDDM) }

1. Initial values of $g_S^0$, $g_{S,\tau}^0$ and $g_D^0$ are guessed, and their values could be taken
zero. We also need to assume the initial values of $\varphi^0_D$ and $\varphi^0_M$ ( to completely decouple the fully mixed dual-permeability system).

2. For $n=1,2,\cdots $, solve the Stokes, mixed microfracture flow and mixed matrix flow equations with Robin-type boundary conditions independently, such that for any $(\mathbf{v}_S,q_S)\in (\mathbf{X}_{S},Q_S)$, $(\mathbf{v}_D,\psi_D)\in (\mathbf{X}_{D},Q_D)$ and $(\mathbf{v}_M,\psi_M)\in (\mathbf{X}_M,Q_M)$, the solutions
$(\mathbf{u}_{S}^{n},p_S^{n})\in (\mathbf{X}_{S}, Q_{S})$, $(\mathbf{u}_{D}^{n},\varphi_D^{n})\in(\mathbf{X}_{D},Q_{D})$ and $(\mathbf{u}_M^n,\varphi_M^n)
\in (\mathbf{X}_M,Q_M)$ can be computed separately from
\begin{eqnarray}
a_{S}(\mathbf{u}^n_{S},\mathbf{v}_{S})-b_S(\mathbf{v}_S, p^n_S)
+\delta_S\langle\mathbf{u}^n_S\cdot\mathbf{n}_S,&\mathbf{v}_{S}&\cdot\mathbf{n}_S\rangle \nonumber\\
&=&(\mathbf{f}_S,\mathbf{v}_S)_{S}-\langle g^{n-1}_S,\mathbf{v}_S\cdot\mathbf{n}_S \rangle
-\sum_{j=1}^{d-1}\langle g^{n-1}_{S,\tau_j},\mathbf{v}_S\cdot \tau_{j} \rangle,\hspace{4mm}\label{decoupled-StokesR1} \\
b_S(\mathbf{u}^n_S, q_S)&=&0,\label{decoupled-StokesR2}\\
a_{D}(\mathbf{u}^n_{D},\mathbf{v}_{D})-b_D(\mathbf{v}_D,\varphi^n_D)
&+&\uwave{\frac{\sigma k_M}{\rho k_D}(\varphi^n_D-\varphi^{n-1}_M, \nabla \cdot \mathbf{v}_D)_D}
+\delta_D\langle\mathbf{u}^n_D\cdot\mathbf{n}_D,\mathbf{v}_D\cdot\mathbf{n}_D\rangle\nonumber\\
&=&\frac{\mu}{\rho k_D}(f_D,\mathrm{div} \hspace{0.5mm} \mathbf{v}_D)_{D}
-\langle g^{n-1}_D,\mathbf{v}_D\cdot\mathbf{n}_D\rangle,\label{decoupled-DarcyR1} \\
b_D(\mathbf{u}^n_D,\psi_D)
+\uwave{\frac{\sigma k_M}{\rho \mu}(\varphi^n_D-\varphi^{n-1}_M, \psi_D)_D}
&=&\frac{1}{\rho}(f_D,\psi_D)_{D},\label{decoupled-DarcyR2}\\
a_{M}(\mathbf{u}^n_{M},\mathbf{v}_{M})-b_M(\mathbf{v}_M,\varphi^n_M)
&+&\uwave{\frac{\sigma }{\rho }(\varphi^n_M-\varphi^{n-1}_D, \nabla \cdot \mathbf{v}_M)_D} = 0,\label{decoupled-DarcyR3} \\
b_M(\mathbf{u}^n_M,\psi_M)
+\uwave{\frac{\sigma k_M}{\rho \mu}(\varphi^n_M-\varphi^{n-1}_D, \psi_M)_D}
&=& 0.\label{decoupled-DarcyR4}
\end{eqnarray}

3. $g_S^{n},\ g_{S,\tau}^{n}$ and $g_D^{n}$ are updated by the following manner:
\begin{eqnarray}
&&g_D^{n} =g_S^{n-1}+(\delta_S+\delta_D)
 {\mathbf{u}}_{S}^n\cdot
\mathbf{n}_{S},\label{decoupled-comp1}\\
&&g_S^{n} =g_D^{n-1}+(\delta_S+\delta_D)
{\mathbf{u}}_{D}^n\cdot \mathbf{n}_{D},\label{decoupled-comp2}\\
&&g_{S,\tau_j}^{n}=-\frac{\nu\alpha\sqrt{d}}{\sqrt{\mathrm{trace}(\mathbf{\prod})}}\mathbf{u}_D^n \cdot \tau_{j} \hspace{6mm} 1\leq j \leq d-1.\label{decoupled-comp3}
\end{eqnarray}

To our best knowledge, this PDDM algorithm is the first completely decoupled scheme for solving the steady-state Stokes-dual-permeability system, which is mainly achieved by the marked terms in (\ref{decoupled-DarcyR1})-(\ref{decoupled-DarcyR4}). It is clearly to see that the existence and uniqueness of above decoupled system solutions $(\mathbf{u}_S^n, p_S^n)$, $(\mathbf{u}_D^n, \varphi_D^n)$ and $(\mathbf{u}_M^n, \varphi_M^n)$ in each iteration follow immediately, because the solutions in this algorithm satisfy the Stokes equation and the mixed Darcy-like equation, respectively. We present the convergence theorem below for this iterative method and demonstrate the convergence of PDDM by applying the elegant energy method.

We first introduce a general convergence lemma for steady-state problems, which is very important for the later convergence analysis of the completely decoupled system.
\begin{lemma}\label{abc}
	Suppose $a_1,a_2,b_1,b_2,c_1,c_2$ are positive constants with $a_2 < a_1,\ b_2 < b_1, c_2 < c_1$. And $A^n,\ B^n, \ C^n \ (n=1,2,\cdots)$ are three different iterative sequence norms, if $a_1 A^n + b_1 B^n + c_1 C^n \leq a_2 A^{n-1} + b_2 B^{n-1} + c_2 C^{n-1}$, we can get
	\begin{eqnarray*}
		&&a_1 A^n + b_1 B^n + c_1 C^n \leq
		\max\Big\{\frac{a_2}{a_1},\frac{b_2}{b_1},\frac{c_2}{c_1} \Big\}^{n-1}\Big(a_2 A^{0} + b_2 B^{0} + c_2 C^{0}\Big).
	\end{eqnarray*}
\end{lemma}
\begin{proof}
	Without loss of generality, we can assume that $\frac{a_2}{a_1}$ is the largest one in $\Big\{ \frac{a_2}{a_1},\frac{b_2}{b_1},\frac{c_2}{c_1} \Big\}$. Using iterative techniques, we can get
	\begin{eqnarray*}
		a_1 A^n + b_1 B^n + c_1 C^n \leq \frac{a_2}{a_1} \Big( a_1 A^{n-1} + b_1 B^{n-1} + c_1 c^{n-1} \Big) 
		+ \Big(b_2-\frac{a_2}{a_1}b_1 \Big) B^{n-1}
		+ \Big(c_2-\frac{a_2}{a_1}c_1 \Big) C^{n-1}.
	\end{eqnarray*}
	Since $\frac{a_2}{a_1}$ is the largest one, we have $\frac{a_2}{a_1}\geq \frac{b_2}{b_1}$ and $\frac{a_2}{a_1}\geq \frac{c_2}{c_1}$, which can induce $b_2-\frac{a_2}{a_1}b_1\leq 0$ and  $c_2-\frac{a_2}{a_1}c_1\leq 0$. Then we can obtain
	\begin{eqnarray*}
		a_1 A^n + b_1 B^n + c_1 C^n &\leq&  \frac{a_2}{a_1} \Big( a_1 A^{n-1} + b_1 B^{n-1} + c_1 C^{n-1} \Big) 
		\leq \frac{a_2}{a_1} \Big( a_2 A^{n-2} + b_2 B^{n-2} + c_2 C^{n-2} \Big)\\ 
		&\leq&\Big( \frac{a_2}{a_1}\Big)^2 \Big( a_1 A^{n-2} + b_1 B^{n-2} + c_1 C^{n-2} \Big) \leq \cdots\\
		&\leq&\Big(\frac{a_2}{a_1} \Big)^{n-1}\Big(a_2 A^{0} + b_2 B^{0} + c_2 C^{0}\Big),
	\end{eqnarray*}
	which proves this lemma.
\end{proof}

\begin{rem}
This lemma is much general but novel, and such lemma or its variants with more different iterative sequence norms (which could be further proved by mathematical induction) can be applied to other steady-state problems. For example, the similar lemma is developed for the ensemble DDM algorithm of the random Stokes-Darcy model \cite{Shi22}. Such lemma plays a similar role with the Gronwall lemma, which is usually used for the unsteady explicit-implicit decoupled algorithm, but fails for the steady-state system.
\end{rem}

\begin{theorem}
Assume that the solution of each iteration (\ref{decoupled-StokesR1})-(\ref{decoupled-DarcyR4}) in the PDDM algorithm is $(\mathbf{u}_S^n,p_S^n;\mathbf{u}_D^n,\varphi_D^n;\mathbf{u}_M^n,\varphi_M^n)$. Let $(\mathbf{u}_S,p_S;\mathbf{u}_D,\varphi_D;\mathbf{u}_M,\varphi_M)$ denotes the solution of the weak formulation (\ref{StokesR1})-(\ref{DarcyR4}). If $\delta_S\leq\delta_D$ and $ g_S^{n},\ g_{S,\tau_j}^{n},\ g_D^{n}$ satisfy the updated compatibility conditions (\ref{decoupled-comp1})-(\ref{decoupled-comp3}), then $(\mathbf{u}_S^n,p_S^n;\mathbf{u}_D^n,\varphi_D^n;\mathbf{u}_M^n,\varphi_M^n)$ converges to $(\mathbf{u}_S,p_S;\mathbf{u}_D,\varphi_D;\mathbf{u}_M,\varphi_M)$.
\end{theorem}

\begin{proof}
Define the following notations for the error functions:
\begin{eqnarray*}
&&\mathbf{e}^n_{S}=\mathbf{u}_{S}-\mathbf{u}^n_{S},\hspace{9mm}\mathbf{e}^n_{D}=\mathbf{u}_{D}-\mathbf{u}^n_{D},
\hspace{13.5mm}\mathbf{e}^n_{M}=\mathbf{u}_{M}-\mathbf{u}^n_{M},\\
&&  {\varepsilon}_S^{n}=p_S-{p}_S^{n},\hspace{10mm}
{\varepsilon}_{D}^{n}=\varphi_D-{\varphi}_D^{n},\hspace{13.5mm}
{\varepsilon}_{M}^{n}=\varphi_M-{\varphi}_M^{n},\\
&&  {\eta}_S^{n}=g_S-{g}^{n}_{S},\hspace{10mm}
{\eta}_{S,\tau_j}^{n}=g_{S,\tau_j}-{g}^{n}_{s,\tau_j}, \hspace{7mm}
{\eta}_D^{n}=g_D-{g}^{n}_{D}.
\end{eqnarray*}
Then for any $(\mathbf{v}_S,q_S;\mathbf{v}_D,\psi_D;\mathbf{v}_M,\psi_M)\in (\mathbf{X}_{S},Q_S;\mathbf{X}_D,Q_D;\mathbf{X}_M,Q_M)$, we can obtain the following error equations by subtracting
(\ref{decoupled-StokesR1})-(\ref{decoupled-DarcyR4}) from
(\ref{StokesR1})-(\ref{DarcyR2}):
\begin{eqnarray}
a_{S}(\mathbf{e}^n_{S},\mathbf{v}_{S})-b_S(\mathbf{v}_S, \varepsilon^n_S)
+\delta_S\langle\mathbf{e}^n_S\cdot\mathbf{n}_S,\mathbf{v}_{S}\cdot\mathbf{n}_S\rangle&&\nonumber\\
=-\langle \eta^{n-1}_S,\mathbf{v}_S\cdot\mathbf{n}_S \rangle
&-&\sum_{j=1}^{d-1}\langle \eta^{n-1}_{S,\tau_j},\mathbf{v}_S\cdot \tau_{j} \rangle,\label{err-1} \\
b_S(\mathbf{e}^n_S, q_S)=0,\hspace{23mm}&&\label{err-2}\\
a_{D}(\mathbf{e}^n_{D},\mathbf{v}_{D})-b_D(\mathbf{v}_D,\varepsilon^n_D)
+\frac{\sigma k_M}{\rho k_D}(\varepsilon^n_D-\varepsilon^{n-1}_M, \nabla \cdot \mathbf{v}_D)_D
&+&\delta_D\langle\mathbf{e}^n_D\cdot\mathbf{n}_D,\mathbf{v}_D\cdot\mathbf{n}_D\rangle\nonumber\\
&=& -\langle \eta^{n-1}_D,\mathbf{v}_D\cdot\mathbf{n}_D\rangle,\label{err-3} \\
b_D(\mathbf{e}^n_D,\psi_D)
+\frac{\sigma k_M}{\rho \mu}(\varepsilon^n_D-\varepsilon^{n-1}_M, \psi_D)_D &=&0,\label{err-4}\\
a_{M}(\mathbf{e}^n_{M},\mathbf{v}_{M})-b_M(\mathbf{v}_M,\varepsilon^n_M)
+\frac{\sigma }{\rho }(\varepsilon^n_M-\varepsilon^{n-1}_D, \nabla \cdot \mathbf{v}_M)_D &=& 0,\label{err-5} \\
b_M(\mathbf{e}^n_M,\psi_M)
+\frac{\sigma k_M}{\rho \mu}(\varepsilon^n_M-\varepsilon^{n-1}_D, \psi_M)_D
&=& 0.\label{err-6}
\end{eqnarray}
Along the interface $\Gamma$, the error functions ${\eta}_D^{n},\ {\eta}_S^{n} $ and ${\eta}_{S,\tau_j}^{n}$ can be updated as follows
\begin{eqnarray}
&&{\eta}_D^{n} ={\eta}_S^{n-1}+(\delta_S+\delta_D)
 {\mathbf{e}}_{S}^n\cdot
\mathbf{n}_{S},\label{errep}\\
&&{\eta}_S^{n} ={\eta}_D^{n-1}+(\delta_S+\delta_D)
{\mathbf{e}}_{D}^n\cdot \mathbf{n}_{D},\label{erref}\\
&&{\eta}_{S,\tau_j}^{n}=-\frac{\nu\alpha\sqrt{d}}{\sqrt{\mathrm{trace}(\mathbf{\prod})}}
\mathbf{e}_D^n \cdot \tau_{j} \hspace{6mm} 1\leq j \leq d-1.\label{erref2}
\end{eqnarray}
Equation (\ref{errep}) can lead to
\begin{eqnarray}\label{etaD2}
||\eta_D^{n}||^2_{\Gamma} =||\eta_S^{n-1}||^2_{\Gamma}+2(\delta_S+\delta_D)\langle\eta_S^{n-1},\mathbf{e}_{S}^{n}\cdot
\mathbf{n}_{S}\rangle+(\delta_S+\delta_D)^2||\mathbf{e}_{S}^{n}\cdot\mathbf{n}_{S}||^2_{\Gamma}.
\end{eqnarray}
Choosing $(\mathbf{v}_{S},q_S)=(\mathbf{e}_{S}^{n},\varepsilon_S^{n})$ in (\ref{err-1})-(\ref{err-2}) and adding the resulting equations together with (\ref{erref2}), we can get
\begin{eqnarray}\label{aS+bS}
a_{S}(\mathbf{e}_{S}^n,\mathbf{e}_{S}^n)+\delta_S||\mathbf{e}_{S}^{n}\cdot\mathbf{n}_{S}||^2_{\Gamma}
=-\langle{\eta}_S^{n-1},\mathbf{e}_S^n\cdot \mathbf{n}_S\rangle+
\frac{\nu\alpha\sqrt{d}}{\sqrt{\mathrm{trace}(\mathbf{\prod})}}
\sum_{j=1}^{d-1}\langle\mathbf{e}_D^{n-1}\cdot \tau_{j},\mathbf{e}_S^n\cdot \tau_{j}\rangle.
\end{eqnarray}
Combining (\ref{etaD2}) and (\ref{aS+bS}), it yields that
\begin{eqnarray}\label{errorD}
||\eta_D^{n}||^2_{\Gamma}
&=& ||\eta_S^{n-1}||^2_{\Gamma}-2(\delta_S+\delta_D)a_S(\mathbf{e}_S^n,\mathbf{e}_S^n)
+(\delta_D^2-\delta_S^2)||\mathbf{e}_{S}^{n}\cdot\mathbf{n}_{S}||^2_{\Gamma}\nonumber\\
&&+2(\delta_S+\delta_D)\frac{\nu\alpha\sqrt{d}}{\sqrt{\mathrm{trace}(\mathbf{\prod})}}
\sum_{j=1}^{d-1}\langle\mathbf{e}_D^{n-1}\cdot \tau_{j},\mathbf{e}_S^n\cdot \tau_{j}\rangle\nonumber\\
&=& ||\eta_S^{n-1}||^2_{\Gamma}-4\nu(\delta_S+\delta_D)||\mathbb{D}(\mathbf{e}_S^n)||_S^2
+(\delta_D^2-\delta_S^2)||\mathbf{e}_{S}^{n}\cdot\mathbf{n}_{S}||^2_{\Gamma}\nonumber\\
&&-2(\delta_S+\delta_D)\frac{\nu\alpha\sqrt{d}}{\sqrt{\mathrm{trace}(\mathbf{\prod})}}
\sum_{j=1}^{d-1}\langle(\mathbf{e}_S^n-\mathbf{e}_D^{n-1})\cdot \tau_{j},\mathbf{e}_S^n\cdot \tau_{j}\rangle.
\end{eqnarray}
Similarly, we can obtain the following equation directly from (\ref{erref}).
\begin{eqnarray}\label{etaS2}
||\eta_S^{n}||^2_{\Gamma}
= ||\eta_D^{n-1}||^2_{\Gamma}+2(\delta_S+\delta_D)\langle\eta_D^{n-1}, \mathbf{e}_D^n\cdot\mathbf{n}_D  \rangle
+(\delta_S+\delta_D)^2||\mathbf{e}_{D}^{n}\cdot\mathbf{n}_{D}||^2_{\Gamma}.
\end{eqnarray}
For the error equations (\ref{err-3})-(\ref{err-4}), we can choose the test function as $(\mathbf{v}_{D},\psi_D)=(\mathbf{e}_{D}^{n},\varepsilon_D^{n})$. Moreover, in order to get a proper analysis, we also need to add the error equations (\ref{err-5})-(\ref{err-6}) with the test function $(\mathbf{v}_{M},\psi_M)=(\mathbf{e}_{M}^{n},\varepsilon_M^{n})$. Then we have
\begin{eqnarray}\label{aD+bD+aM+bM}
&&a_{D}(\mathbf{e}_{D}^n,\mathbf{e}_{D}^n)+\frac{\sigma k_M}{\rho k_D}(\varepsilon_D^n-\varepsilon_M^{n-1},\nabla\cdot\mathbf{e}_D^n)_D
+\frac{\sigma k_M}{\rho \mu}(\varepsilon_D^n-\varepsilon_M^{n-1},\varepsilon_D^n)_D
+\delta_D||\mathbf{e}_{D}^{n}\cdot\mathbf{n}_{D}||^2_{\Gamma}\nonumber\\
&&\hspace{4mm}+a_{M}(\mathbf{e}_{M}^n,\mathbf{e}_{M}^n)
+\frac{\sigma }{\rho } (\varepsilon_M^n-\varepsilon_D^{n-1},\nabla\cdot\mathbf{e}_M^n)_D
+\frac{\sigma k_M}{\rho \mu}(\varepsilon_M^n-\varepsilon_D^{n-1},\varepsilon_M^n)_D=-\langle{\eta}_D^{n-1},\mathbf{e}_D^n\cdot \mathbf{n}_D\rangle.\hspace{6mm}
\end{eqnarray}
Combining (\ref{etaS2}) and (\ref{aD+bD+aM+bM}), we can obtain
\begin{eqnarray}\label{errorS}
||\eta_S^{n}||^2_{\Gamma}
= ||\eta_D^{n-1}||^2_{\Gamma}&-&2(\delta_S+\delta_D)a_D(\mathbf{e}_D^n,\mathbf{e}_D^n)
-2(\delta_S+\delta_D)\frac{\sigma k_M}{\rho k_D} (\varepsilon_D^n-\varepsilon_M^{n-1},\nabla\cdot\mathbf{e}_D^n)_D
\nonumber\\
&-&2(\delta_S+\delta_D)\frac{\sigma k_M}{\rho \mu}(\varepsilon_D^n-\varepsilon_M^{n-1},\varepsilon_D^n)_D
+(\delta_S^2-\delta_D^2)||\mathbf{e}_{D}^{n}\cdot\mathbf{n}_{D}||^2_{\Gamma}\nonumber\\
&-&2(\delta_S+\delta_D)a_{M}(\mathbf{e}_{M}^n,\mathbf{e}_{M}^n)
-2(\delta_S+\delta_D)\frac{\sigma }{\rho } (\varepsilon_M^n-\varepsilon_D^{n-1},\nabla\cdot\mathbf{e}_M^n)_D
\nonumber\\
&-&2(\delta_S+\delta_D)\frac{\sigma k_M}{\rho \mu}(\varepsilon_M^n-\varepsilon_D^{n-1},\varepsilon_M^n)_D.
\end{eqnarray}

Since the convergence analysis for cases $\delta_S=\delta_D$ and $\delta_S<\delta_D$ are different, we will treat them separately.

\textbf{Case 1: $\delta_S=\delta_D=\delta$.} In this case, we can simplify the equations (\ref{errorD}) and (\ref{errorS}) to arrive at
\begin{eqnarray}
\label{errorD-1}
||\eta_D^{n}||^2_{\Gamma}
&=& ||\eta_S^{n-1}||^2_{\Gamma}-8\nu\delta||\mathbb{D}(\mathbf{e}_S^n)||_S^2
-4\delta\frac{\nu\alpha\sqrt{d}}{\sqrt{\mathrm{trace}(\mathbf{\prod})}}
\sum_{j=1}^{d-1}\langle(\mathbf{e}_S^n-\mathbf{e}_D^{n-1})\cdot \tau_{j},\mathbf{e}_S^n\cdot \tau_{j}\rangle,\\
\label{errorS-1}
||\eta_S^{n}||^2_{\Gamma}
&=& ||\eta_D^{n-1}||^2_{\Gamma}-4\delta a_D(\mathbf{e}_D^n,\mathbf{e}_D^n)
-4\delta\frac{\sigma k_M}{\rho k_D} (\varepsilon_D^n-\varepsilon_M^{n-1},\nabla\cdot\mathbf{e}_D^n)_D
-4\delta\frac{\sigma k_M}{\rho \mu}(\varepsilon_D^n-\varepsilon_M^{n-1},\varepsilon_D^n)_D
\nonumber\\
&&-4\delta a_{M}(\mathbf{e}_{M}^n,\mathbf{e}_{M}^n)
-4\delta\frac{\sigma }{\rho } (\varepsilon_M^n-\varepsilon_D^{n-1},\nabla\cdot\mathbf{e}_M^n)_D
-4\delta\frac{\sigma k_M}{\rho \mu} (\varepsilon_M^n-\varepsilon_D^{n-1},\varepsilon_M^n)_D.
\end{eqnarray}
Adding equations (\ref{errorD-1})-(\ref{errorS-1}) together, then summing over $n$ from $n=1$ to $N$, we can get
\begin{eqnarray}\label{etaS1+etaD1}
||\eta_S^{N}||^2_{\Gamma}+||\eta_D^{N}||^2_{\Gamma}
&=&||\eta_S^{0}||^2_{\Gamma}+||\eta_D^{0}||^2_{\Gamma}
-4\delta \sum_{n=1}^{N}\Big[2\nu||\mathbb{D}(\mathbf{e}_S^n)||_S^2
+ a_D(\mathbf{e}_D^n,\mathbf{e}_D^n) + a_M(\mathbf{e}_M^n,\mathbf{e}_M^n) \nonumber\\
&&+\frac{\nu\alpha\sqrt{d}}{\sqrt{\mathrm{trace}(\mathbf{\prod})}}
\sum_{j=1}^{d-1}\langle(\mathbf{e}_S^n-\mathbf{e}_D^{n-1})\cdot \tau_{j},\mathbf{e}_S^n\cdot \tau_{j}\rangle\nonumber\\
&&+\frac{\sigma k_M}{\rho k_D} (\varepsilon_D^n-\varepsilon_M^{n-1},\nabla\cdot\mathbf{e}_D^n)_D
+\frac{\sigma k_M}{\rho \mu}(\varepsilon_D^n-\varepsilon_M^{n-1},\varepsilon_D^n)_D\nonumber\\
&&+\frac{\sigma }{\rho } (\varepsilon_M^n-\varepsilon_D^{n-1},\nabla\cdot\mathbf{e}_M^n)_D
+\frac{\sigma k_M}{\rho \mu}(\varepsilon_M^n-\varepsilon_D^{n-1},\varepsilon_M^n)_D \Big].\hspace{8mm}
\end{eqnarray}
Recalling some trace inequalities \cite{trace}, 
there exist constants $C_{\mathrm{tr}}$, $C'_{\mathrm{tr}}$, $C''_{\mathrm{tr}}$, which only depend on the domain $\Omega_D$ or $\Omega_S$, such that for any $\mathbf{v}_D\in \mathbf{X}_D$ or $\mathbf{v}_S\in \mathbf{X}_S$,
\begin{eqnarray}\label{trace}
	||\mathbf{v}_D||_{H^{-\frac{1}{2}}(\Gamma)}\leq C_{\mathrm{tr}} ||\mathbf{v}_D||_{\mathrm{div}}, \hspace{3mm} ||\mathbf{v}_S||_{H^{\frac{1}{2}}(\Gamma)}\leq C'_{\mathrm{tr}} ||\mathbf{v}_S||_{1}, \hspace{3mm} 
	||\mathbf{v}_S||_{\Gamma}\leq C''_{\mathrm{tr}} ||\mathbf{v}_S||_{S}^{\frac{1}{2}}||\mathbf{v}_S||_{1}^{\frac{1}{2}}.\hspace{3mm}
\end{eqnarray}
Thanks to Cauchy-Schwarz inequality, trace inequality (\ref{trace}), Poincar\'e inequality, Korn's inequality, and Young's inequality, there exist positive constants $C_1,\ C_2$, 
such that
\begin{eqnarray}
&&\sum_{n=1}^N\Big[2\nu||\mathbb{D}(\mathbf{e}_S^n)||_S^2
+ a_D(\mathbf{e}_D^n,\mathbf{e}_D^n) + a_M(\mathbf{e}_M^n,\mathbf{e}_M^n)
+\frac{\nu\alpha\sqrt{d}}{\sqrt{\mathrm{trace}(\mathbf{\prod})}}
\sum_{j=1}^{d-1}\langle(\mathbf{e}_S^n-\mathbf{e}_D^{n-1})\cdot \tau_{j},\mathbf{e}_S^n\cdot \tau_{j}\rangle\Big]\nonumber\\
&&\geq\sum_{n=1}^N\Big[2\nu||\mathbb{D}(\mathbf{e}_S^n)||_S^2
+\frac{\mu}{\rho k_D} ||\mathbf{e}_D||^2_{\mathrm{div}}
+\frac{\mu }{\rho k_M} ||\mathbf{e}_M||^2_{\mathrm{div}}\nonumber\\
&&\hspace{12mm}+\frac{\nu\alpha\sqrt{d}}{\sqrt{\mathrm{trace}(\mathbf{\prod})}}
\sum_{j=1}^{d-1}||\mathbf{e}_S^n\cdot \tau_{j}||_{\Gamma}^2
-\frac{\nu\alpha\sqrt{d}}{\sqrt{\mathrm{trace}(\mathbf{\prod})}}
\sum_{j=1}^{d-1}||\mathbf{e}_D^{n-1}\cdot \tau_{j}||_{H^{-\frac{1}{2}}(\Gamma)}
||\mathbf{e}_S^n\cdot \tau_{j}||_{H^{\frac{1}{2}}(\Gamma)}\Big]\nonumber\\
&&\geq \sum_{n=1}^N\Big[C_1\nu||\mathbf{e}_S^n||_{1}^2
+\frac{\mu}{\rho k_D} ||\mathbf{e}_D^n||_{\mathrm{div}}^2+\frac{\mu}{\rho k_M}||\mathbf{e}_M^n||_{\mathrm{div}}^2\nonumber\\
&&\hspace{12mm}-\frac{(C'_{\mathrm{tr}})^2\nu \alpha d}{2 C_1 \mathrm{trace}(\mathbf{\prod})}\sum_{j=1}^{d-1}||\mathbf{e}_D^{n-1}\cdot \tau_{j}||^2_{H^{-\frac{1}{2}}(\Gamma)}
-\frac{C_1 \nu}{2 (C'_{\mathrm{tr}})^2}\sum_{j=1}^{d-1}||\mathbf{e}_S^n\cdot \tau_{j}||^2_{H^{\frac{1}{2}}(\Gamma)}
\Big]\nonumber\\
&&\geq \sum_{n=1}^N\Big[C_1\nu||\mathbf{e}_S^n||_{1}^2
+\frac{\mu}{\rho k_D} ||\mathbf{e}_D^n||_{\mathrm{div}}^2+\frac{\mu}{\rho k_M} ||\mathbf{e}_M^n||_{\mathrm{div}}^2
-\frac{C_1\nu}{2}||\mathbf{e}_S^n||_{1}^2
-\frac{C_2^2\nu\alpha^2d}{2C_1\mathrm{trace}(\mathbf{\prod})}
||\mathbf{e}_D^{n}||_{\mathrm{div}}^2\Big]\nonumber\\
&&\hspace{4mm}-\frac{(C'_{\mathrm{tr}})^2\nu \alpha d}{2 C_1 \mathrm{trace}(\mathbf{\prod})}\sum_{j=1}^{d-1}||\mathbf{e}_D^{0}\cdot \tau_{j}||^2_{H^{-\frac{1}{2}}(\Gamma)}\nonumber\\
&&\geq\frac{C_1\nu}{2}\sum_{n=1}^N||\mathbf{e}_S^n||_{1}^2
+\Big(\frac{\mu}{\rho k_D} -\frac{C_2^2\nu\alpha^2d}{2C_1\mathrm{trace}(\mathbf{\prod})}\Big)
\sum_{n=1}^N||\mathbf{e}_D^n||_{\mathrm{div}}^2
+\frac{\mu}{\rho k_M}\sum_{n=1}^N||\mathbf{e}_M^n||_{\mathrm{div}}^2
-\frac{(C'_{\mathrm{tr}})^2 }{2 C_1 \nu} \sum_{j=1}^{d-1}||\eta^0_{S,\tau_j}||^2_{H^{-\frac{1}{2}}(\Gamma)}.\hspace{9mm}
\label{firstterm}
\end{eqnarray}
Applying the equality condition $(a-b,a)=\frac{1}{2}(||a||^2-||b||^2+||a-b||^2)$, Cauchy-Schwarz inequality, and Young's inequality, we obtain
\begin{eqnarray}
&&\frac{\sigma k_M}{\rho\mu}\sum_{n=1}^N\Big[
(\varepsilon_D^n-\varepsilon_M^{n-1},\varepsilon_D^n)_D
+(\varepsilon_M^n-\varepsilon_D^{n-1},\varepsilon_M^n)_D
+\frac{\mu}{k_D}(\varepsilon_D^n-\varepsilon_M^{n-1},\nabla\cdot\mathbf{e}_D^n)_D
+\frac{\mu}{k_M}(\varepsilon_M^n-\varepsilon_D^{n-1},\nabla\cdot\mathbf{e}_M^n)_D\Big]\nonumber\\
&&= \frac{\sigma k_M}{2\rho\mu}\sum_{n=2}^N\Big[||\varepsilon_D^n||_D^2-||\varepsilon_M^{n-1}||_D^2
+||\varepsilon_D^n-\varepsilon_M^{n-1}||_D^2
+||\varepsilon_M^n||_D^2-||\varepsilon_D^{n-1}||_D^2
+||\varepsilon_M^n-\varepsilon_D^{n-1}||_D^2\Big]\nonumber\\
&&\hspace{4mm}+\frac{\sigma k_M}{\rho\mu}\sum_{n=1}^N
\Big[\frac{\mu}{k_D}(\varepsilon_D^n-\varepsilon_M^{n-1},\nabla\cdot\mathbf{e}_D^n)_D
+\frac{\mu}{k_M}(\varepsilon_M^n-\varepsilon_D^{n-1},\nabla\cdot\mathbf{e}_M^n)_D\Big]\nonumber\\
&&\geq\frac{\sigma k_M}{2\rho\mu}\Big[||\varepsilon_D^N||_D^2+||\varepsilon_M^{N}||_D^2
-(||\varepsilon_D^0||_D^2+||\varepsilon_M^{0}||_D^2)\Big]
+\frac{\sigma k_M}{2\rho\mu}\sum_{n=1}^N\Big[||\varepsilon_D^n-\varepsilon_M^{n-1}||_D^2
+||\varepsilon_M^n-\varepsilon_D^{n-1}||_D^2\Big]\nonumber\\
&&\hspace{4mm}-\frac{\sigma k_M}{\rho k_D}\sum_{n=1}^N
||\varepsilon_D^n-\varepsilon_M^{n-1}||_D||\nabla\cdot\mathbf{e}_D^n||_D
-\frac{\sigma }{\rho}\sum_{n=1}^N
||\varepsilon_M^n-\varepsilon_D^{n-1}||_D||\nabla\cdot\mathbf{e}_M^n||_D\nonumber\\
&&\geq\frac{\sigma k_M}{2\rho\mu}\Big[||\varepsilon_D^N||_D^2+||\varepsilon_M^{N}||_D^2
-(||\varepsilon_D^0||_D^2+||\varepsilon_M^{0}||_D^2)\Big]
+\frac{\sigma k_M}{2\rho\mu}\sum_{n=1}^N\Big[||\varepsilon_D^n-\varepsilon_M^{n-1}||_D^2
+||\varepsilon_M^n-\varepsilon_D^{n-1}||_D^2\Big]\nonumber\\
&&\hspace{4mm}-\frac{\sigma k_M}{2\rho\mu}\sum_{n=1}^N||\varepsilon_D^n-\varepsilon_M^{n-1}||_D^2
-\frac{\sigma \mu k_M}{2\rho k_D^2} \sum_{n=1}^N||\nabla\cdot\mathbf{e}_D^n||_D^2 -\frac{\sigma k_M}{2\rho\mu} \sum_{n=1}^N||\varepsilon_M^n-\varepsilon_D^{n-1}||_D^2
-\frac{\sigma \mu}{2\rho k_M} \sum_{n=1}^N||\nabla\cdot\mathbf{e}_M^n||_D^2\nonumber\\
&&\geq\frac{\sigma k_M}{2\rho\mu}\Big[||\varepsilon_D^N||_D^2+||\varepsilon_M^{N}||_D^2
-(||\varepsilon_D^0||_D^2+||\varepsilon_M^{0}||_D^2)\Big]
-\frac{\sigma \mu k_M}{2\rho k_D^2} \sum_{n=1}^N||\mathbf{e}_D^n||_{\mathrm{div}}^2
-\frac{\sigma \mu}{2\rho k_M} \sum_{n=1}^N||\mathbf{e}_M^n||_{\mathrm{div}}^2.
\label{secondterm}
\end{eqnarray}
To this end, substituting (\ref{firstterm}) and (\ref{secondterm}) into (\ref{etaS1+etaD1}), we have the desired estimate as follows
\begin{eqnarray*}\label{converge0}
0&\leq& ||\eta_S^{N}||^2_{\Gamma}+||\eta_D^{N}||^2_{\Gamma}\nonumber\\
&\leq& ||\eta_S^{0}||^2_{\Gamma}+||\eta_D^{0}||^2_{\Gamma}
+\frac{2\delta\sigma k_M}{\rho\mu}(||\varepsilon_D^0||_D^2+||\varepsilon_M^{0}||_D^2)
+ \frac{2\delta(C'_{\mathrm{tr}})^2 }{ C_1 \nu} \sum_{j=1}^{d-1}||\eta^0_{S,\tau_j}||^2_{H^{-\frac{1}{2}}(\Gamma)} \nonumber\\
&&-4\delta \sum_{n=1}^N \Big[ \frac{\nu C_1}{2}||\mathbf{e}_S^n||_{1}^2
+\Big(\frac{2\mu k_D-\sigma \mu k_M}{2 \rho k_D^2}
-\frac{C_2^2\nu\alpha^2d}{2C_1\mathrm{trace}(\mathbf{\prod})}\Big)||\mathbf{e}_D^n||_{\mathrm{div}}^2
+\frac{\mu (2-\sigma) }{2\rho k_M} ||\mathbf{e}_M^n||_{\mathrm{div}}^2 \Big].\hspace{10mm}
\end{eqnarray*}
The intrinsic permeability $k_D$ and $k_M$ in microfractures and matrix are generally selected as a value less than or equal to $1$ and normally $k_D \gg k_M$. 
Following \cite{Cao10}, we also assume the parameter $\alpha$ is smaller, so that $\frac{2\mu k_D-\sigma \mu k_M}{2 \rho k_D^2}
-\frac{C_2^2\nu\alpha^2d}{2C_1\mathrm{trace}(\mathbf{\prod})}>0$ and $2-\sigma>0$.
For any positive integer $N$, we can finally get
\begin{eqnarray*}
&&4\delta \sum_{n=1}^N \Big[ \frac{\nu C_1}{2}||\mathbf{e}_S^n||_{1}^2
+\Big(\frac{2\mu k_D-\sigma \mu k_M}{2 \rho k_D^2}
-\frac{C_2^2\nu\alpha^2d}{2C_1\mathrm{trace}(\mathbf{\prod})}\Big)||\mathbf{e}_D^n||_{\mathrm{div}}^2
+\frac{\mu (2-\sigma) }{2\rho k_M} ||\mathbf{e}_M^n||_{\mathrm{div}}^2 \Big]\\
&&\hspace{3mm} \leq   ||\eta_S^{0}||^2_{\Gamma}+||\eta_D^{0}||^2_{\Gamma}
+\frac{2\delta\sigma k_M}{\rho\mu}(||\varepsilon_D^0||_D^2+||\varepsilon_M^{0}||_D^2)
+ \frac{2\delta(C'_{\mathrm{tr}})^2 }{ C_1 \nu} \sum_{j=1}^{d-1}||\eta^0_{S,\tau_j}||^2_{H^{-\frac{1}{2}}(\Gamma)}.
\end{eqnarray*}
Therefore, by the convergence theorem, $\mathbf{e}_{S}^{n}$, $\mathbf{e}_{D}^{n}$ and $\mathbf{e}_M^n$ tend to be zero in $H^1(\Omega_S)^d$, $H(\mathrm{div};\Omega_D)$ and $H(\mathrm{div};\Omega_D)$ respectively.

The convergence of series $||\varepsilon_{S}^{n}||_S$ will be proved next. By utilizing the similar result as \cite{Sun2021}, for given $\varepsilon_{S}^{n}\in Q_S$, there exists
$\mathbf{v}^{\varepsilon}_{S}\in \mathbf{X}_S \cap H_0^1(\Omega_S)^d$ and a positive constant $C_{I}$, such that the following results can be obtained,
\begin{eqnarray*}
\nabla\cdot \mathbf{v}^{\varepsilon}_{S}=\varepsilon_{S}^{n}\ \
\mathrm{in} \ \Omega_S, \hspace{10mm} || \mathbf{v}^{\varepsilon}_{S}||_{1}\leq
C_{I}||\varepsilon_{S}^{n}||_S.
\end{eqnarray*}
Since $\mathbf{v}^{\varepsilon}_{S}\in (H^1_0(\Omega_S)^d)$, it is clearly to see  $\mathbf{v}^{\varepsilon}_{S}=0$ on the boundary $ \partial \Omega_S$. Besides this, we can check easily that $\mathbf{v}^{\varepsilon}_{S} \in \mathbf{X}_S$ and ${\mathbf{v}}_{S}^{\varepsilon}\cdot {\mathbf{n}}_{S}=0, \ {\mathbf{v}}_{S}^{\varepsilon}\cdot {\mathbf{\tau}}_{j}=0$ on the interface $\Gamma$. Through the above premise, an inequality can be received as follows,
\begin{eqnarray*}
b_{S}(\mathbf{v}_{S}^{\varepsilon},\varepsilon_S^{n})
=||\varepsilon_S^{n}||_S^2\geq
1/C_{I}||\mathbf{v}_{S}^{\varepsilon}||_{1}||\varepsilon_S^{n}||_S.
\end{eqnarray*}
Let the test function $\mathbf{v}_S={\mathbf{v}}_{S}^{\varepsilon}$
in the error equation (\ref{err-1}). Since ${\mathbf{v}}_{S}^{\varepsilon}\cdot {\mathbf{n}}_{S}=0$
and ${\mathbf{v}}_{S}^{\varepsilon}\cdot {\mathbf{\tau}}_{j}=0$ on the interface $\Gamma$,
the following equation can be inferred,
\begin{eqnarray*}\label{err-1r}
a_{S}(\mathbf{e}_{S}^{n},\mathbf{v}_{S}^{\varepsilon})-b_{S}(\mathbf{v}%
_{S}^{\varepsilon},\varepsilon_S^{n})=0.
\end{eqnarray*}
Finally, thanks to Cauchy-Schwarz inequality, there exists a positive constant $C_3$ such that
\begin{eqnarray}\label{S_P}
 ||\varepsilon_S^{n}||_S \leq C_{I}\frac{b_{S}(\mathbf{v}%
_{S}^{\varepsilon},\varepsilon_S^{n})}{||\mathbf{v}_{S}^{\varepsilon}||_{1}}
=C_{I}\frac{a_{S}(\mathbf{e}_{S}^{n},\mathbf{v}_{S}^{\varepsilon})}
{||\mathbf{v}_{S}^{\varepsilon}||_{1}}
\leq C_3\nu||\mathbf{e}_{S}^{n}||_{1}.
\end{eqnarray}
This implies the convergence of $||\varepsilon_S^{n}||_S$, which means $\varepsilon_S^{n}$ tends to be zero in $L^2_0(\Omega_S)$.

For the convergence of series $||\varepsilon_D^{n}||_D$ and $||\varepsilon_M^{n}||_D$, we prove the convergence of $||\varepsilon_D^{n}-\varepsilon_M^{n-1}||_D$ and $||\varepsilon_M^{n}-\varepsilon_D^{n-1}||_D$  firstly.
Due to $\varepsilon_D^{n}-\varepsilon_M^{n-1} \in L_0^2(\Omega_D)$ and $\varepsilon_M^{n}-\varepsilon_D^{n-1} \in L_0^2(\Omega_D)$, we choose $\psi_D = \varepsilon_D^{n}-\varepsilon_M^{n-1} $ in (\ref{err-4}) and $\psi_M = \varepsilon_M^{n}-\varepsilon_D^{n-1} $ in (\ref{err-6}) to get:
\begin{eqnarray*}
	&&\frac{\sigma k_M}{\rho \mu} ||\varepsilon_D^{n}-\varepsilon_M^{n-1} ||_D^2 =-b_D(\mathbf{e}_D^n,\varepsilon_D^{n}-\varepsilon_M^{n-1})
	\leq \frac{1}{\rho} ||\mathbf{e}_D^n||_{\mathrm{div}} ||\varepsilon_D^{n}-\varepsilon_M^{n-1}||_D,\\
	&&\frac{\sigma k_M}{\rho \mu} ||\varepsilon_M^{n}-\varepsilon_D^{n-1} ||_D^2 =-b_M(\mathbf{e}_M^n,\varepsilon_M^{n}-\varepsilon_D^{n-1})
	\leq \frac{1}{\rho} ||\mathbf{e}_M^n||_{\mathrm{div}} ||\varepsilon_M^{n}-\varepsilon_D^{n-1}||_D.
\end{eqnarray*}
Then, we have
\begin{eqnarray}\label{D-M}
	&& ||\varepsilon_D^{n}-\varepsilon_M^{n-1} ||_D
	\leq \frac{\mu}{\sigma k_M} ||\mathbf{e}_D^n||_{\mathrm{div}},\\
	&& ||\varepsilon_M^{n}-\varepsilon_D^{n-1} ||_D
	\leq \frac{\mu}{\sigma k_M} ||\mathbf{e}_M^n||_{\mathrm{div}}.\label{M-D}
\end{eqnarray}
Hence, we have obtained the convergence of $||\varepsilon_D^{n}-\varepsilon_M^{n-1}||_D$ and $||\varepsilon_M^{n}-\varepsilon_D^{n-1}||_D$.

For $\varepsilon_{D}^{n}\in Q_D$ and $\varepsilon_{M}^{n}\in Q_M$, there exists $\mathbf{{v}}_D^{\varepsilon} \in \mathbf{X}_D \cap H_0^1(\Omega_D)^d$, $\mathbf{{v}}_M^{\varepsilon} \in \mathbf{X}_M \cap H_0^1(\Omega_D)^d$ and positive constants $C_{II}, C_{III}$ satisfying
\begin{eqnarray*}
&&\nabla\cdot \mathbf{{v}}^{\varepsilon}_{D}=\varepsilon_{D}^{n}\ \
\mathrm{in} \ \Omega_D, \hspace{7.5mm}
\mathbf{{v}}^{\varepsilon}_{D}=0 \ \ \mathrm{on} \ \partial\Omega_D, \hspace{6.5mm}
|| \mathbf{{v}}^{\varepsilon}_{D}||_{1}\leq
C_{II}||\varepsilon_{D}^{n}||_D,\\
&&\nabla\cdot \mathbf{{v}}^{\varepsilon}_{M}=\varepsilon_{M}^{n}\ \
\mathrm{in} \ \Omega_D, \hspace{6mm}
\mathbf{{v}}^{\varepsilon}_{M}=0 \ \ \mathrm{on} \ \partial\Omega_D, \hspace{6mm}
|| \mathbf{{v}}^{\varepsilon}_{M}||_{1}\leq
C_{III}||\varepsilon_{M}^{n}||_D.
\end{eqnarray*}
Furthermore, we can directly get
\begin{eqnarray*}
&&b_{D}(\mathbf{{v}}_{D}^{\varepsilon},\varepsilon_D^{n})
=\frac{1}{\rho}||\varepsilon_D^{n}||_D^2
\geq \frac{1}{\rho C_{II}}||\mathbf{{v}}_{D}^{\varepsilon}||_{1}||\varepsilon_D^{n}||_D
\geq \frac{1}{\rho C_{II}} ||\mathbf{{v}}_{D}^{\varepsilon}||_{\mathrm{div}}||\varepsilon_D^{n}||_D,\\
&&b_{M}(\mathbf{{v}}_{M}^{\varepsilon},\varepsilon_M^{n})
=\frac{1}{\rho}||\varepsilon_M^{n}||_D^2
\geq \frac{1}{\rho C_{III}} ||\mathbf{{v}}_{M}^{\varepsilon}||_{1}||\varepsilon_M^{n}||_D
\geq \frac{1}{\rho C_{III}} ||\mathbf{{v}}_{M}^{\varepsilon}||_{\mathrm{div}}||\varepsilon_M^{n}||_D.
\end{eqnarray*}
By selecting the test functions $\mathbf{v}_D={\mathbf{{v}}}_{D}^{\varepsilon}$ and $\mathbf{v}_M={\mathbf{{v}}}_{M}^{\varepsilon}$ for the error equations (\ref{err-3}) and (\ref{err-5}) respectively, and using the property of the test functions to be $0$ on the boundary $\partial \Omega_D$, we can get the following equations:
\begin{eqnarray*}
&&a_{D}(\mathbf{e}_{D}^{n},\mathbf{{v}}_{D}^{\varepsilon})
-b_{D}(\mathbf{{v}}_{D}^{\varepsilon},\varepsilon_D^{n})
+\frac{\sigma k_M}{\rho k_D}(\varepsilon^n_D-\varepsilon^{n-1}_M,
\nabla \cdot \mathbf{{v}}_D^{\varepsilon})_D=0,\\
&&a_{M}(\mathbf{e}_{M}^{n},\mathbf{{v}}_{M}^{\varepsilon})
-b_{M}(\mathbf{{v}}_{M}^{\varepsilon},\varepsilon_M^{n})
+\frac{\sigma }{\rho }(\varepsilon^n_M-\varepsilon^{n-1}_D,
\nabla \cdot \mathbf{{v}}_M^{\varepsilon})_D=0.
\end{eqnarray*}
The convergence of series $||\varepsilon_D^{n}||_D$ and $||\varepsilon_M^{n}||_D$ can be derived by the Cauchy-Schwarz inequality, which is shown as follows:
\begin{eqnarray}\label{Dcontrol}
||\varepsilon_D^{n}||_D
&\leq& \rho C_{II} \frac{b_{D}(\mathbf{{v}}_{D}^{\varepsilon},\varepsilon_D^{n})}
{||\mathbf{{v}}_{D}^{\varepsilon}||_{\mathrm{div}}}
=\frac{\rho C_{II} a_{D}(\mathbf{e}_{D}^{n},\mathbf{\tilde{v}}_{D}^{\varepsilon})
+\frac{\sigma k_M C_{II}}{\mu K_D} (\varepsilon^n_D-\varepsilon^{n-1}_M,
\nabla \cdot \mathbf{{v}}_D^{\varepsilon})_D}
{||\mathbf{{v}}_{D}^{\varepsilon}||_{\mathrm{div}}} \nonumber\\
&\leq& \frac{ \mu C_{II}}{k_D}  ||\mathbf{e}_{D}^{n}||_{\mathrm{div}}
+\frac{\sigma k_M C_{II}}{ k_D} ||\varepsilon^n_D-\varepsilon^{n-1}_M||_D
\leq \frac{ 2 \mu C_{II}}{k_D}  ||\mathbf{e}_{D}^{n}||_{\mathrm{div}},\\
||\varepsilon_M^{n}||_D
&\leq& \rho C_{III} \frac{b_{M}(\mathbf{{v}}_{M}^{\varepsilon},\varepsilon_M^{n})}
{||\mathbf{{v}}_{M}^{\varepsilon}||_{\mathrm{div}}}
=\frac{\rho C_{III} a_{M}(\mathbf{e}_{M}^{n},\mathbf{{v}}_{M}^{\varepsilon})
+\sigma C_{III} (\varepsilon^n_M-\varepsilon^{n-1}_D,
\nabla \cdot \mathbf{{v}}_M^{\varepsilon})_D}
{||\mathbf{{v}}_{M}^{\varepsilon}||_{\mathrm{div}}}\nonumber\\
&\leq& \frac{\mu}{k_M} C_{III} ||\mathbf{e}_{M}^{n}||_{\mathrm{div}}
+\sigma C_{III} ||\varepsilon^n_M-\varepsilon^{n-1}_D||_D
\leq   \frac{2\mu C_{III}}{k_M} ||\mathbf{e}_{M}^{n}||_{\mathrm{div}},\label{Mcontrol}
\end{eqnarray}
Therefore, $\varepsilon_D^{n}$ and $\varepsilon_M^{n}$ tend to be zero in $L_0^2(\Omega_D)$.

Next, we will prove that $\eta_{S,\tau_j}^n$ ( $ j\in[1,d-1]$ ) tends to be zero in $H^{-\frac{1}{2}}(\Gamma)$. Utilizing the equation (\ref{erref2}) and trace inequality (\ref{trace}) with a positive constant $C_4$, we can yield
\begin{eqnarray}\label{eta_Stau}
||{\eta}_{S,\tau_j}^{n}||_{H^{-\frac{1}{2}}(\Gamma)}
\leq\frac{\nu\alpha\sqrt{d}}{\sqrt{\mathrm{trace}(\mathbf{\prod})}}
||\mathbf{e}_D^{n} \cdot \tau_{j}||_{H^{-\frac{1}{2}}(\Gamma)}
\leq\frac{C_4\nu\alpha\sqrt{d}}{\sqrt{\mathrm{trace}(\mathbf{\prod})}}
||\mathbf{e}_D^n ||_{\mathrm{div}}.
\end{eqnarray}
So $||{\eta}_{S,\tau_j}^{n}||_{H^{-\frac{1}{2}}(\Gamma)}$ is convergent.

In order to demonstrate the convergence of $\eta_D^n$ and $\eta_S^n$, we can basically refer to the analysis of \cite{Sun2021} (Page 9). A similar auxiliary function $\lambda=\delta \mathbf{e}_D^{n}\cdot \mathbf{n}_D + \eta_D^{n-1}$ on $\Gamma$ can be constructed. After the same treatment in \cite{Sun2021}, we can get the following inequality
$\beta_{\Gamma} ||\lambda||_{{0},\Gamma}
\leq \frac{<\delta \mathbf e^n_D\cdot \mathbf n_D+\eta_D^{n-1}, \mathbf w_D\cdot \mathbf
n_D>}{||\mathbf w_D||_{\mathrm{div}}},$
where $\beta_{\Gamma}<1$ is a positive constant, $\mathbf w_D$ is a test function for equation (\ref{err-3}), and the detail of $\mathbf w_D$ can be found in \cite{Sun2021}. With above results, we can further obtain
\begin{eqnarray*}
\beta_{\Gamma} ||\lambda||_{{0},\Gamma}
&\leq& \frac{<\delta \mathbf e^n_D\cdot \mathbf n_D+\eta_D^{n-1}, \mathbf w_D\cdot \mathbf n_D>}
{||\mathbf w_D||_{\mathrm{div}}}\\
&\leq&\frac{\left| a_{D}(\mathbf{e}_{D}^{n},\mathbf{w}_{D})
-b_{D}(\mathbf{w}_{D},\varepsilon_D^{n})
+\frac{\sigma k_M}{\rho k_D}(\varepsilon^n_D-\varepsilon^{n-1}_M,
\nabla \cdot \mathbf{w}_D)_D \right|}
{||\mathbf{w}_D||_{\mathrm{div}}}\\
&\leq&\frac{\mu}{\rho k_D}||\mathbf{e}_{D}^{n}||_{\mathrm{div}}
+\frac{1}{\rho}||\varepsilon_D^{n}||_D
+\frac{\sigma k_M}{\rho k_D} ||\varepsilon^n_D-\varepsilon^{n-1}_M||_D.
\end{eqnarray*}
Combining with the triangle inequality and trace inequality, it follows that
\begin{eqnarray*}
||\eta_D^{n-1}||_{H^{-\frac{1}{2}}(\Gamma)} &\leq&  ||\lambda||_{H^{-\frac{1}{2}}(\Gamma)}+ \delta ||\mathbf e_D^{n}\cdot \mathbf{n}_D||_{H^{-\frac{1}{2}}(\Gamma)}\\
&\leq& C_4 ||\lambda||_{\Gamma}+ \delta ||\mathbf e_D^{n}\cdot \mathbf
n_D||_{H^{-\frac{1}{2}}(\partial \Omega_D)}\\
&\leq& \frac{C_5\tilde{C}_D+C_6\rho\delta}{\rho}||\mathbf{e}_{D}^{n}||_{\mathrm{div}}
+\frac{C_5}{\rho}||\varepsilon_D^{n}||_D
+\frac{C_5 \sigma k_M}{\rho k_D} ||\varepsilon^n_D-\varepsilon^{n-1}_M||_D,
\end{eqnarray*}
with some positive constants $C_5, C_6$. Due to ${\eta}_S^{n} ={\eta}_D^{n-1}+2\delta
{\mathbf{e}}_{D}^n\cdot \mathbf{n}_{D}$, we finally have
\begin{eqnarray*}
||\eta_S^{n}||_{H^{-\frac{1}{2}}(\Gamma)}&\leq&  ||\eta_D^{n-1}||_{H^{-\frac{1}{2}}(\Gamma)}+ 2\delta ||\mathbf e_D^{k}\cdot \mathbf n_D||_{H^{-\frac{1}{2}}(\Gamma)}\\
&\leq& ||\eta_D^{n-1}||_{H^{-\frac{1}{2}}(\Gamma)}+ 2\delta||\mathbf e_D^{n}\cdot \mathbf
n_D||_{H^{-\frac{1}{2}}(\partial \Omega_D)}\\
&\leq& \frac{C_5\tilde{C}_D+3C_6\rho\delta}{\rho}||\mathbf{e}_{D}^{n}||_{\mathrm{div}}
+\frac{C_5}{\rho}||\varepsilon_D^{n}||_D
+\frac{C_5 \sigma k_M}{\rho k_D} ||\varepsilon^n_D-\varepsilon^{n-1}_M||_D.
\end{eqnarray*}
Therefore, the convergence of $\eta_D^n$ and $\eta_S^n$ in $H^{-1/2}(\Gamma)$ is guaranteed.

\textbf{Case 2: $\delta_S<\delta_D$.} In this case, we need to recall the equations (\ref{errorD}) and (\ref{errorS}) for further analysis. 
Using equation (\ref{errorS}) to replace the term $||\eta_S^{n-1}||_{\Gamma}^2$ in equation (\ref{errorD}), we can have
\begin{eqnarray*}
||\eta_D^{n}||^2_{\Gamma}
&=&||\eta_D^{n-2}||^2_{\Gamma}+(\delta_D^2-\delta_S^2)||\mathbf{e}_{S}^{n}\cdot\mathbf{n}_{S}||^2_{\Gamma}
+(\delta_S^2-\delta_D^2)||\mathbf{e}_{D}^{n-1}\cdot\mathbf{n}_{D}||^2_{\Gamma}\\
&&-2(\delta_S+\delta_D)\Big[
2\nu||\mathbb{D}(\mathbf{e}_S^n)||_S^2+a_D(\mathbf{e}_D^{n-1},\mathbf{e}_D^{n-1})
+a_{M}(\mathbf{e}_{M}^{n-1},\mathbf{e}_{M}^{n-1})\\
&&\hspace{22mm}+\frac{\nu\alpha\sqrt{d}}{\sqrt{\mathrm{trace}(\mathbf{\prod})}}
\sum_{j=1}^{d-1}\langle(\mathbf{e}_S^n-\mathbf{e}_D^{n-1})\cdot \tau_{j},\mathbf{e}_S^n\cdot \tau_{j}\rangle\\
&&\hspace{22mm}+\frac{\sigma k_M}{\rho k_D} (\varepsilon_D^{n-1}-\varepsilon_M^{n-2},\nabla\cdot\mathbf{e}_D^{n-1})_D
+\frac{\sigma k_M}{\rho \mu} (\varepsilon_D^{n-1}-\varepsilon_M^{n-2},\varepsilon_D^{n-1})_D\\
&&\hspace{22mm}+\frac{\sigma }{\rho } (\varepsilon_M^{n-1}-\varepsilon_D^{n-2},\nabla\cdot\mathbf{e}_M^{n-1})_D
+\frac{\sigma k_M}{\rho \mu}(\varepsilon_M^{n-1}-\varepsilon_D^{n-2},\varepsilon_M^{n-1})_D\Big]\\
&=&||\eta_D^{n-2}||^2_{\Gamma}
+(\delta_S^2-\delta_D^2)||\mathbf{e}_{D}^{n-1}\cdot\mathbf{n}_{D}||^2_{\Gamma}\\
&&-2(\delta_S+\delta_D)\Big[
2\nu||\mathbb{D}(\mathbf{e}_S^n)||_S^2+a_D(\mathbf{e}_D^{n-1},\mathbf{e}_D^{n-1})+a_{M}(\mathbf{e}_{M}^{n-1},\mathbf{e}_{M}^{n-1})\\
&&\hspace{21mm}+\frac{\nu\alpha\sqrt{d}}{\sqrt{\mathrm{trace}(\mathbf{\prod})}}
\sum_{j=1}^{d-1}\langle(\mathbf{e}_S^n-\mathbf{e}_D^{n-1})\cdot \tau_{j},\mathbf{e}_S^n\cdot \tau_{j}\rangle
-\frac{\delta_D-\delta_S}{2}||\mathbf{e}_S^n\cdot\mathbf{n}_S||_{\Gamma}^2\\
&&\hspace{21mm}+\frac{\sigma k_M}{\rho k_D} (\varepsilon_D^{n-1}-\varepsilon_M^{n-2},\nabla\cdot\mathbf{e}_D^{n-1})_D
+\frac{\sigma k_M}{\rho \mu} (\varepsilon_D^{n-1}-\varepsilon_M^{n-2},\varepsilon_D^{n-1})_D\\
&&\hspace{21mm}+\frac{\sigma }{\rho} (\varepsilon_M^{n-1}-\varepsilon_D^{n-2},\nabla\cdot\mathbf{e}_M^{n-1})_D
+\frac{\sigma k_M}{\rho \mu} (\varepsilon_M^{n-1}-\varepsilon_D^{n-2},\varepsilon_M^{n-1})_D\Big].
\end{eqnarray*}
By utilizing Korn's inequality, Cauchy-Schwarz inequality, Young's inequality and trace inequality, the following derivation process is carried out,
\begin{eqnarray*}
&&2\nu||\mathbb{D}(\mathbf{e}_S^n)||_S^2
+ a_D(\mathbf{e}_D^{n-1},\mathbf{e}_D^{n-1}) 
+\frac{\nu\alpha\sqrt{d}}{\sqrt{\mathrm{trace}(\mathbf{\prod})}}
\sum_{j=1}^{d-1}\langle(\mathbf{e}_S^n-\mathbf{e}_D^{n-1})\cdot \tau_{j},\mathbf{e}_S^n\cdot \tau_{j}\rangle
-\frac{\delta_D-\delta_S}{2}||\mathbf{e}_S^n\cdot\mathbf{n}_S||_{\Gamma}^2\\
&&\hspace{4mm}\geq C_1\nu||\mathbf{e}_S^n||_{1}^2
+\frac{\mu }{\rho k_D} ||\mathbf{e}_D^{n-1}||_{\mathrm{div}}^2
-\Big(\frac{C_1\nu}{2}||\mathbf{e}_S^n||_{1}^2
+\frac{C_2^2\nu\alpha^2d}{2C_1\mathrm{trace}(\mathbf{\prod})}
||\mathbf{e}_D^{n-1}||_{\mathrm{div}}^2\Big)
-\frac{C_7(\delta_D-\delta_S)}{2}||\mathbf{e}_S^n||_{1}^2\\
&&\hspace{4mm}\geq \frac{C_1\nu-C_7(\delta_D-\delta_S)}{2}||\mathbf{e}_S^n||_{1}^2
+\Big(\frac{\mu}{\rho k_D} -\frac{C_2^2\nu\alpha^2d}{2C_1\mathrm{trace}(\mathbf{\prod})}\Big)
||\mathbf{e}_D^{n-1}||_{\mathrm{div}}^2,
\end{eqnarray*}
with a positive constant $C_7$, and other constants have been defined before. Then similar to the derivation of inequality (\ref{secondterm}), we can apply  $(a-b,a)=\frac{1}{2}(||a||^2-||b||^2+||a-b||^2)$, and thanks to the inequality conditions (\ref{Dcontrol})-(\ref{Mcontrol}), Cauchy-Schwarz inequality, and Young's inequality, we can continue to get
\begin{eqnarray*}
&&\frac{\sigma k_M}{\rho \mu}\Big[(\varepsilon_D^{n-1}-\varepsilon_M^{n-2},\varepsilon_D^{n-1})_D
+(\varepsilon_D^{n-1}-\varepsilon_M^{n-2},\varepsilon_D^{n-1})_D \Big]\\
&&\hspace{12mm}+\frac{\sigma k_M}{\rho k_D} (\varepsilon_D^{n-1}-\varepsilon_M^{n-2},\nabla\cdot\mathbf{e}_D^{n-1})_D
+\frac{\sigma }{\rho} (\varepsilon_M^{n-1}-\varepsilon_D^{n-2},\nabla\cdot\mathbf{e}_M^{n-1})_D\\
&&\hspace{8mm} \geq \frac{\sigma k_M}{2 \rho \mu} \Big[
||\varepsilon_D^{n-1}||_D^2-||\varepsilon_M^{n-2}||_D^2+||\varepsilon_D^{n-1}-\varepsilon_M^{n-2}||_D^2
+||\varepsilon_M^{n-1}||_D^2-||\varepsilon_D^{n-2}||_D^2+||\varepsilon_M^{n-1}-\varepsilon_D^{n-2}||_D^2
\Big]\\
&&\hspace{12mm} -\frac{\sigma k_M}{2 \rho \mu} ||\varepsilon_D^{n-1}-\varepsilon_M^{n-2}||_D^2 
- \frac{\mu \sigma  k_M}{2 \rho k_D^2}||\mathbf{e}_D^{n-1}||_{\mathrm{div}}^2
-\frac{\sigma k_M}{2 \rho \mu} ||\varepsilon_M^{n-1}-\varepsilon_D^{n-2}||_D^2
-\frac{\mu \sigma  }{2 \rho k_M}||\mathbf{e}_M^{n-1}||_{\mathrm{div}}^2\\
&&\hspace{8mm} \geq -\frac{\sigma k_M}{2 \rho \mu} ||\varepsilon_D^{n-2}||_D^2 
-\frac{\sigma k_M}{2 \rho \mu} ||\varepsilon_M^{n-2}||_D^2
-\frac{\mu \sigma  k_M}{2 \rho k_D^2}||\mathbf{e}_D^{n-1}||_{\mathrm{div}}^2
-\frac{\mu \sigma  }{2 \rho k_M}||\mathbf{e}_M^{n-1}||_{\mathrm{div}}^2\\
&&\hspace{8mm} \geq -\frac{\mu \sigma  k_M}{2 \rho k_D^2}||\mathbf{e}_D^{n-1}||_{\mathrm{div}}^2
-\frac{2 \mu \sigma  k_M C_{II}^2}{ \rho k_D^2}||\mathbf{e}_D^{n-2}||_{\mathrm{div}}^2
-\frac{\mu \sigma  }{2 \rho k_M} ||\mathbf{e}_M^{n-1}||_{\mathrm{div}}^2
-\frac{2 \mu \sigma C_{III}^2 }{ \rho k_M}||\mathbf{e}_M^{n-2}||_{\mathrm{div}}^2.
\end{eqnarray*}
In summary, we can yield
\begin{eqnarray}\label{laststep}
||\eta_D^{n}||^2_{\Gamma}
&\leq&||\eta_D^{n-2}||^2_{\Gamma}
+(\delta_S^2-\delta_D^2)||\mathbf{e}_{D}^{n-1}\cdot\mathbf{n}_{D}||^2_{\Gamma}-2(\delta_S+\delta_D)\Big[\frac{C_1\nu-C_7(\delta_D-\delta_S)}{2}||\mathbf{e}_S^n||_{1}^2\nonumber\\
&&+\Big(\frac{2\mu k_D-\mu \sigma  k_M}{2\rho k_D^2} -\frac{C_2^2\nu\alpha^2d}{2C_1\mathrm{trace}(\mathbf{\prod})}\Big)
||\mathbf{e}_D^{n-1}||_{\mathrm{div}}^2-\frac{2 \mu \sigma  k_M C_{II}^2}{ \rho k_D^2}||\mathbf{e}_D^{n-2}||_{\mathrm{div}}^2\nonumber\\
&&+\frac{ \mu ( 2-\sigma )  }{2 \rho k_M} ||\mathbf{e}_M^{n-1}||_{\mathrm{div}}^2-\frac{2 \mu \sigma C_{III}^2 }{ \rho k_M}||\mathbf{e}_M^{n-2}||_{\mathrm{div}}^2\Big].
\end{eqnarray}
By adapting the same technique utilized in \cite{Sun2021, GR86}, we can obtain a test function $\mathbf{\tilde{v}}_D\in \mathbf{X}_D$, which satisfies $\mathbf{\tilde{v}}\cdot\mathbf{n}_D|_{\Gamma}=\eta_D^{n-2}$ and $||\mathbf{\tilde{v}}_D||_{\mathrm{div}}\leq||\eta_D^{n-2}||_{\Gamma}$.
Substituting $\mathbf{\tilde{v}}_D$ into (\ref{err-3}), and using the Cauchy-Schwarz inequality, Young's inequality and (\ref{Dcontrol}), we have
\begin{eqnarray}\label{contral_eta}
||\eta_{D}^{n-2}||_{\Gamma}^{2}&=&
-a_{D}(\mathbf{e}^{n-1}_{D},\mathbf{\tilde{v}}_{D})+b_D(\mathbf{\tilde{v}}_D,\varepsilon^{n-1}_D)
+\frac{\sigma k_M}{\rho k_D} (\varepsilon^{n-1}_D-\varepsilon^{n-2}_M, \nabla \cdot \mathbf{\tilde{v}}_D)_D
+\delta_D\langle\mathbf{e}^{n-1}_D\cdot\mathbf{n}_D,\eta_D^{n-2}\rangle\nonumber\\
&\leq&\frac{\mu}{\rho k_D} ||\mathbf{e}_{D}^{n-1}||_{{\mathrm{div}}}||\eta_D^{n-2}||_{\Gamma}
+ \frac{ 2 \mu C_{II}}{\rho k_D}  ||\mathbf{e}_{D}^{n-1}||_{\mathrm{div}}
||\eta_D^{n-2}||_{\Gamma}
+\frac{\mu}{\rho k_D} ||\mathbf{e}_D^{n-1}||_{\mathrm{div}}||\eta_D^{n-2}||_{\Gamma}\nonumber\\
&&\hspace{36.6mm}
+\delta_D||\mathbf{e}_{D}^{n-1}\cdot\mathbf{n}_D||_{\Gamma} ||\eta_{D}^{n-2}||_{\Gamma}\nonumber\\
&\leq&
\frac{\mu (\delta_S+\delta_D)}{\rho k_D} ||\mathbf{e}_D^{n-1}||_{\mathrm{div}}^2
+(\delta_D^2-\delta_S^2)||\mathbf{e}_{D}^{n-1}\cdot\mathbf{n}_{D}||^2_{\Gamma}
+\theta(\delta_S,\delta_D) ||\eta_{D}^{n-2}||_{\Gamma}^2,
\end{eqnarray}
where 
\begin{eqnarray}
\theta(\delta_S,\delta_D):=\frac{4(C_{II}+1)^2(\delta_D-\delta_S)-\rho k_D \delta_D}{4\rho k_D (\delta_D^2-\delta_S^2)}.
\end{eqnarray}
Combining (\ref{laststep}) and (\ref{contral_eta}), we obtain
\begin{eqnarray*}
	||\eta_D^{n}||^2_{\Gamma}
	&\leq&\theta(\delta_S,\delta_D)||\eta_D^{n-2}||^2_{\Gamma}
	-2(\delta_S+\delta_D)\Big[\frac{C_1\nu-C_7(\delta_D-\delta_S)}{2}||\mathbf{e}_S^n||_{1}^2\nonumber\\
	&&+\Big(\frac{\mu k_D-\mu \sigma  k_M}{2\rho k_D^2} -\frac{C_2^2\nu\alpha^2d}{2C_1\mathrm{trace}(\mathbf{\prod})}\Big)
	||\mathbf{e}_D^{n-1}||_{\mathrm{div}}^2-\frac{2 \mu \sigma  k_M C_{II}^2}{ \rho k_D^2}||\mathbf{e}_D^{n-2}||_{\mathrm{div}}^2\nonumber\\
	&&+\frac{ \mu ( 2-\sigma )  }{2 \rho k_M} ||\mathbf{e}_M^{n-1}||_{\mathrm{div}}^2-\frac{2 \mu \sigma C_{III}^2 }{ \rho k_M} ||\mathbf{e}_M^{n-2}||_{\mathrm{div}}^2\Big].
\end{eqnarray*}
Then suppose the parameter $\alpha$ is small enough and
\begin{eqnarray}\label{final-Conditions}
	\delta_D-\delta_S < \frac{\mu C_1}{C_7},\hspace{6mm}
	\theta(\delta_S,\delta_D)<1,\hspace{6mm}
	k_D \gg k_M, \hspace{6mm} 
	\sigma <  \frac{2}{4C_{III}^2+1}.
\end{eqnarray}
Refer to \cite{Sun2021} (Remark 3.4), we can choose the appropriate parameters $\delta_S $ and $\delta_D$ to get the condition $\theta(\delta_S,\delta_D)<1$. Then, we can have
\begin{eqnarray*}
	&&\frac{C_1\nu-C_7(\delta_D-\delta_S)}{2}>0,\hspace{10mm} \frac{ \mu ( 2-\sigma )  }{2 \rho k_M} > \frac{2 \mu \sigma C_{III}^2 }{ \rho k_M}>0, \\
	&&\frac{\mu k_D-\mu \sigma  k_M}{2\rho k_D^2} -\frac{C_2^2\nu\alpha^2d}{2C_1\mathrm{trace}(\mathbf{\prod})} >\frac{2 \mu \sigma  k_M C_{II}^2}{ \rho k_D^2}>0,	
\end{eqnarray*}
so that
\begin{eqnarray*}
	&&2(\delta_S+\delta_D)\Big[\Big(\frac{\mu k_D-\mu \sigma  k_M}{2\rho k_D^2} -\frac{C_2^2\nu\alpha^2d}{2C_1\mathrm{trace}(\mathbf{\prod})}\Big)
	||\mathbf{e}_D^{n-1}||_{\mathrm{div}}^2
	+ \frac{ \mu ( 2-\sigma )  }{2 \rho k_M} ||\mathbf{e}_M^{n-1}||_{\mathrm{div}}^2\Big]
	+||\eta_D^{n}||^2_{\Gamma}\\
	&&\hspace{16mm}\leq\theta(\delta_S,\delta_D)||\eta_D^{n-2}||^2_{\Gamma}
	+2(\delta_S+\delta_D)\Big[ \frac{2 \mu \sigma  k_M C_{II}^2}{ \rho k_D^2}||\mathbf{e}_D^{n-2}||_{\mathrm{div}}^2
	+\frac{2 \mu \sigma C_{III}^2 }{ \rho k_M} ||\mathbf{e}_M^{n-2}||_{\mathrm{div}}^2\Big].
\end{eqnarray*}
After a small technical transformation, adding the term $\sqrt{\theta(\delta_S,\delta_D)}||\eta_D^{n-1}||^2_{\Gamma}$ to the left and right of the above inequality, it can be obtained
\begin{eqnarray*}
	&&2(\delta_S+\delta_D)\Big[\Big(\frac{\mu k_D-\mu \sigma  k_M}{2\rho k_D^2} -\frac{C_2^2\nu\alpha^2d}{2C_1\mathrm{trace}(\mathbf{\prod})}\Big)
	||\mathbf{e}_D^{n-1}||_{\mathrm{div}}^2
	+ \frac{ \mu ( 2-\sigma )  }{2 \rho k_M} ||\mathbf{e}_M^{n-1}||_{\mathrm{div}}^2\Big]\\
	&&\hspace{20mm}+\Big( ||\eta_D^{n}||^2_{\Gamma}+ \sqrt{\theta(\delta_S,\delta_D)}||\eta_D^{n-1}||^2_{\Gamma} \Big)\\
	&&\hspace{16mm}\leq
	2(\delta_S+\delta_D)\Big[ \frac{2 \mu \sigma  k_M C_{II}^2}{ \rho k_D^2}||\mathbf{e}_D^{n-2}||_{\mathrm{div}}^2
	+\frac{2 \mu \sigma C_{III}^2 }{ \rho k_M} ||\mathbf{e}_M^{n-2}||_{\mathrm{div}}^2\Big]\\
	&&\hspace{20mm}+\sqrt{\theta(\delta_S,\delta_D)}\Big(||\eta_D^{n-1}||^2_{\Gamma}+ \sqrt{\theta(\delta_S,\delta_D)}||\eta_D^{n-2}||^2_{\Gamma} \Big).
\end{eqnarray*}
	By utilizing Lemma \ref{abc}, we get
	\begin{eqnarray*}
		&&2(\delta_S+\delta_D)\Big[\Big(\frac{\mu k_D-\mu \sigma  k_M}{2\rho k_D^2} -\frac{C_2^2\nu\alpha^2d}{2C_1\mathrm{trace}(\mathbf{\prod})}\Big)
		||\mathbf{e}_D^{N-1}||_{\mathrm{div}}^2
		+ \frac{ \mu ( 2-\sigma )  }{2 \rho k_M} ||\mathbf{e}_M^{N-1}||_{\mathrm{div}}^2\Big]\\
		&&\hspace{20mm}+\Big( ||\eta_D^{N}||^2_{\Gamma}+ \sqrt{\theta(\delta_S,\delta_D)}||\eta_D^{N-1}||^2_{\Gamma} \Big)\\
		&&\hspace{16mm}\leq \max\Big\{ \frac{4 \sigma  k_M C_{II}^2}{ k_D- \sigma  k_M 
			-\frac{C_2^2\nu\alpha^2d \rho k_D^2}{\mu C_1\mathrm{trace}(\mathbf{\prod})} }, \frac{ 4 \sigma C_{III}^2}{2-\sigma }, \sqrt{\theta(\delta_S,\delta_D)} \Big\}^{N-2}\\
		&&\hspace{27mm}\Big[2(\delta_S+\delta_D) \frac{2 \mu \sigma  k_M C_{II}^2}{ \rho k_D^2} ||\mathbf{e}_D^{0}||_{\mathrm{div}}^2
		+2(\delta_S+\delta_D)\frac{2 \mu \sigma C_{III}^2 }{ \rho k_M} ||\mathbf{e}_M^{0}||_{\mathrm{div}}^2\\
		&&\hspace{32mm}+\sqrt{\theta(\delta_S,\delta_D)}\Big(||\eta_D^{1}||^2_{\Gamma}+ \sqrt{\theta(\delta_S,\delta_D)}||\eta_D^{0}||^2_{\Gamma} \Big)\Big],
	\end{eqnarray*}
which means $||\mathbf{e}_D^{n}||_{\mathrm{div}}^2$, $||\mathbf{e}_M^{n}||_{\mathrm{div}}^2$ and $||\eta_D^{n}||^2_{\Gamma}$ are geometric convergence.

Note that we add the error equations (\ref{err-5})-(\ref{err-6}) with the test function $(\mathbf{v}_{M},\psi_M)=(\mathbf{e}_{M}^{n},\varepsilon_M^{n})$ to get the equation (\ref{errorS}), in fact, these additional items add up to zero. So we can delete them to simplify equation  (\ref{errorS}), which is shown as follows
\begin{eqnarray}\label{errorS-C2}
	||\eta_S^{n}||^2_{\Gamma}
	&=& ||\eta_D^{n-1}||^2_{\Gamma}-2(\delta_S+\delta_D)a_D(\mathbf{e}_D^n,\mathbf{e}_D^n)
	-2(\delta_S+\delta_D)\frac{\sigma k_M}{\rho k_D} (\varepsilon_D^n-\varepsilon_M^{n-1},\nabla\cdot\mathbf{e}_D^n)_D
	\nonumber\\
	&&-2(\delta_S+\delta_D)\frac{\sigma k_M}{\rho \mu}(\varepsilon_D^n-\varepsilon_M^{n-1},\varepsilon_D^n)_D
	+(\delta_S^2-\delta_D^2)||\mathbf{e}_{D}^{n}\cdot\mathbf{n}_{D}||^2_{\Gamma}.
\end{eqnarray}
The equation (\ref{errorS-C2}) is further analyzed by utilizing (\ref{Mcontrol}), we can get:
\begin{eqnarray*}
	||\eta_S^{n}||^2_{\Gamma}
	&\leq& ||\eta_D^{n-1}||^2_{\Gamma}-2(\delta_S+\delta_D)a_D(\mathbf{e}_D^n,\mathbf{e}_D^n)
	-2(\delta_S+\delta_D)\frac{\sigma k_M}{\rho k_D} (\varepsilon_D^n-\varepsilon_M^{n-1},\nabla\cdot\mathbf{e}_D^n)_D\\
	&&-2(\delta_S+\delta_D)\frac{\sigma k_M}{\rho \mu}(\varepsilon_D^n-\varepsilon_M^{n-1},\varepsilon_D^n)_D\\
	&\leq& ||\eta_D^{n-1}||^2_{\Gamma} + 2(\delta_S+\delta_D) \Big[ 
	-\frac{\mu}{\rho k_D} ||\mathbf{e}_D^n||_{\mathrm{div}}^2 + \frac{\sigma k_M}{2 \rho \mu} ||\varepsilon_D^n-\varepsilon_M^{n-1}||_D^2 
	+ \frac{\mu \sigma k_M}{2 \rho k_D^2} ||\mathbf{e}_D^n||_{\mathrm{div}}^2\\
	&&\hspace{37mm}-\frac{\sigma k_M}{2 \rho \mu} \Big( ||\varepsilon_D^n||_D^2-||\varepsilon_M^{n-1}||_D^2+||\varepsilon_D^n-\varepsilon_M^{n-1}||_D^2 \Big)
	\Big]\\
	&\leq& ||\eta_D^{n-1}||^2_{\Gamma} + 2(\delta_S+\delta_D) \Big[
	\Big( \frac{\mu \sigma k_M-2\mu k_D}{2 \rho k_D^2} \Big) ||\mathbf{e}_D^n||_{\mathrm{div}}^2
	+\frac{\sigma k_M}{2 \rho \mu}||\varepsilon_M^{n-1}||_D^2
	\Big] \\
		&\leq& ||\eta_D^{n-1}||^2_{\Gamma} + 2(\delta_S+\delta_D) \Big[
	 \frac{\mu \sigma k_M-2\mu k_D}{2 \rho k_D^2} ||\mathbf{e}_D^n||_{\mathrm{div}}^2
	+\frac{2 \mu \sigma C_{III}^2}{2 \rho k_M}||\mathbf{e}_M^{n-1}||_{\mathrm{div}}^2
	\Big],
\end{eqnarray*}
where $\frac{\mu \sigma k_M-2\mu k_D}{2 \rho k_D^2}<0$ under the conditions (\ref{final-Conditions}), which means $||\eta_S^{n}||^2_{\Gamma}$ is geometric convergence.

By combining the above geometric convergence conclusion and the Korn's inequality with the equation (\ref{aS+bS}), we obtain
\begin{eqnarray*}
	C_1 \nu ||\mathbf{e}_S^n||_1^2 &\leq& a_{S}(\mathbf{e}_{S}^n,\mathbf{e}_{S}^n)+\delta_S||\mathbf{e}_{S}^{n}\cdot\mathbf{n}_{S}||^2_{\Gamma}
	=-\langle{\eta}_S^{n-1},\mathbf{e}_S^n\cdot \mathbf{n}_S\rangle+
	\frac{\nu\alpha\sqrt{d}}{\sqrt{\mathrm{trace}(\mathbf{\prod})}}
	\sum_{j=1}^{d-1}\langle\mathbf{e}_D^{n-1}\cdot \tau_{j},\mathbf{e}_S^n\cdot \tau_{j}\rangle\\
	&\leq& ||{\eta}_S^{n-1}||_{\Gamma}||\mathbf{e}_S^n||_1
	+\frac{\nu\alpha\sqrt{d}}{\sqrt{\mathrm{trace}(\mathbf{\prod})}}
	||\mathbf{e}_D^{n-1}||_{\mathrm{div}}||\mathbf{e}_S^n||_1,
\end{eqnarray*}
so that
\begin{eqnarray*}
	||\mathbf{e}_S^n||_1 &\leq& \frac{1}{C_1 \nu} ||{\eta}_S^{n-1}||_{\Gamma}
	+\frac{\nu\alpha\sqrt{d}}{C_1 \nu \sqrt{\mathrm{trace}(\mathbf{\prod})}}
	||\mathbf{e}_D^{n-1}||_{\mathrm{div}},
\end{eqnarray*}
which implies the geometric convergence of $||\mathbf{e}_S^n||_1$.

Moreover, combining with the results (\ref{S_P}), (\ref{Dcontrol}),
(\ref{Mcontrol}) and (\ref{eta_Stau}), we can summarize the following geometric convergence result:
	\begin{eqnarray}
	&&||\mathbf{e}_S^N||_1^2+||\varepsilon_S^{N}||_S^2+||\mathbf{e}_D^{N-1}||_{\mathrm{div}}^2 + ||\varepsilon_D^{N-1}||_D^2
	+ ||\mathbf{e}_M^{N-1}||_{\mathrm{div}}^2 + ||\varepsilon_M^{N-1}||_D^2\nonumber\\
	&&\hspace{25.5mm} + ||\eta_S^{N}||^2_{\Gamma} + ||\eta_D^{N}||^2_{\Gamma} +\sqrt{\theta(\delta_S,\delta_D)}||\eta_D^{N-1}||^2_{\Gamma} + \sum_{j=1}^{d-1}||{\eta}_{S,\tau_j}^{N-1}||_{H^{-\frac{1}{2}}(\Gamma)}^2\nonumber\\
	&&\hspace{16mm}\leq C^*\max\Big\{ \frac{4 \sigma  k_M C_{II}^2}{ k_D- \sigma  k_M 
		-\frac{C_2^2\nu\alpha^2d \rho k_D^2}{\mu C_1\mathrm{trace}(\mathbf{\prod})} }, \frac{ 4 \sigma C_{III}^2}{2-\sigma }, \sqrt{\theta(\delta_S,\delta_D)} \Big\}^{N-2}\nonumber\\
	&&\hspace{33mm}\Big[2(\delta_S+\delta_D) \frac{2 \mu \sigma  k_M C_{II}^2}{ \rho k_D^2} ||\mathbf{e}_D^{0}||_{\mathrm{div}}^2
	+2(\delta_S+\delta_D)\frac{2 \mu \sigma C_{III}^2 }{ \rho k_M} ||\mathbf{e}_M^{0}||_{\mathrm{div}}^2\nonumber\\
	&&\hspace{42mm}+\sqrt{\theta(\delta_S,\delta_D)}\Big(||\eta_D^{1}||^2_{\Gamma}+ \sqrt{\theta(\delta_S,\delta_D)}||\eta_D^{0}||^2_{\Gamma} \Big)\Big],
\end{eqnarray}
for a given positive constant $C^*$.
\end{proof}


\section{Numerical Experiments}
In this section, three numerical experiments are presented to illustrate the accuracy and efficiency of the proposed parallel domain decomposition algorithm for the fully-mixed Stokes-dual-permeability fluid flow model with BJ interface conditions. As for the finite element spaces to all experiments, we can choose the well-known MINI (P1b-P1) elements for the Stokes part and Brezzi-Douglas-Marini (BDM1-P0) elements for both microfracture and matrix parts. The stopping criterion for the iterative process of PDDM algorithm is usually selected as: 
\begin{eqnarray*}
	RE_{\mathtt{stop}}:=\Bigl(\frac{||\mathbf{u}^{n}_{S,h}-\mathbf{u}_{S,h}^{n-1}||_S^2}{||(\mathbf{u}^{n}_{S,h}+\mathbf{u}_{S,h}^{n-1})/2||_S^2} +\frac{||\mathbf{u}^{n}_{D,h}-\mathbf{u}_{D,h}^{n-1}||_D^2}{||(\mathbf{u}^{n}_{D,h}+\mathbf{u}_{D,h}^{n-1})/2||_D^2}+\frac{||\mathbf{u}^{n}_{M,h}-\mathbf{u}_{M,h}^{n-1}||_D^2}{||(\mathbf{u}^{n}_{M,h}+\mathbf{u}_{M,h}^{n-1})/2||_D^2} \Bigr)^{1/2}\leq 10^{-6}.
\end{eqnarray*}
Here, the subscript $h$ denotes that solutions are derived by finite element methods and $RE_{\mathtt{stop}}$ is definded as the relative approximation error of the stopping criterion. 

In the first numerical experiment, we test a smooth problem with exact solutions to verify the convergence and check the feasibility of the PDDM algorithm. We present the flow speed and  streamlines on a horizontal cased-hole completion wellbore with a vertical production wellbore in the third example, which is more efficiently in petroleum engineering application. Above numerical tests are implemented by the open software FreeFEM++ \cite{F18}. In the second test, a coupling 3D shallow water system with a complicated dual-permeability region is simulated to show the conservation of mass on the interface. Due to the limitation of 3D finite elements in the software FreeFEM++, another open-source FEniCS \cite{FEniCS} is selected for further implementation.

\subsection{Smooth problem with exact solutions}
The smooth problem with exact solutions is adopted from \cite{HouJY16} to support the convergence analysis and check the optimal error orders of the PDDM algorithm. The dual-porosity region $\Omega_D=[0.0,1.0]\times [0.0,0.75]$ and the fluid region $\Omega_S=[0.0,1.0]\times [-0.25,0.0]$  with the interface $\Gamma=[0.0,1.0]\times\{0\}$ are considered. In order to get the external body force and source term, we should choose the exact solutions as follows:
\begin{eqnarray*}
	\begin{aligned}
		\mathbf{u}_D&= -\frac{k_D}{\mu} \nabla \varphi_D, \hspace{12mm}	\varphi_{D} = (2-\pi\sin(\pi x))(\cos(\pi(1-y))-y),\\
		\mathbf{u}_M&= -\frac{k_M}{\mu} \nabla \varphi_M, \hspace{11mm}	\varphi_{M} =\sin(xy^2-y^3),\\
		\mathbf{u}_{S} &=\Big[x^2y^2+\exp(-y), \ -\frac{2}{3}xy^3+(2-\pi\sin(\pi x)) \Big]^{\mathtt{T}}, \hspace{4mm} p_{S} = (\pi \sin(\pi x)-2)\cos(2\pi y),
	\end{aligned}
\end{eqnarray*}
which can be easily verified that the BJ interface conditions are satisfied. For better comparison with \cite{HouJY16} and computational convenience, we select the physical parameters  $\nu, \sigma, \mu, \rho, \alpha$ as $1.0$ and the intrinsic permeabilities as $k_D=1.0, k_M=0.01$.

We solve this smooth problem by the PDDM algorithm on a uniform triangular mesh with the mesh size $h$. From the theoretical analysis, we can know that the proposed algorithm is not parameter free, the convergence is closely dependent on the selection of Robin parameters $\delta_S$ and $\delta_D$. So, in Fig. 5.1, we plot the variation trend of $RE_{\mathtt{stop}}$ with the increase of iteration steps $n$ for different choices of $\delta_S$ and $\delta_D$, while $h=\frac{1}{64}$. It is clear to observe that $RE_{\mathtt{stop}}$ gradually tends to zero while $\delta_S \leq \delta_D$. Moreover, while the intrinsic permeabilities are selected as $k_D=1.0, k_M=0.01$, the proposed algorithm is divergent in the case $\delta_S > \delta_D$.
\begin{figure}[htbp] 
	\label{REstop}
	\centering
	\includegraphics[width=65mm,height=50mm]{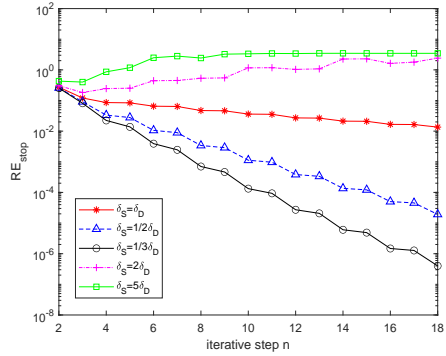}
	\caption{The variation trend of the relative approximation error of the stopping criterion $RE_{\mathtt{stop}}$ while $k_D=1.0, k_M=0.01$ and $h=\frac{1}{64}$. }
\end{figure}

To check the feasibility of PDDM, the relative numerical errors of the velocities and pressures in Stokes, microfracture and matrix are provided in Table \ref{OrderDDM}. One can see clearly that the PDDM algorithm has the optimal orders.
Especially, the numerical errors of velocities in the microfracture and matrix is second-order in $L^2$-norm, which demonstrates that the mixed finite element methods (PDDM) are superior to classical Galerkin methods \cite{HouJY16, HouJY21}. Moreover, the iteration step $n$ increases with mesh refinement while $\delta_S=\delta_D$. However, for the case $\delta_S<\delta_D$, the convergence rate of PDDM is $h$-independent, which supports the theoretical analysis. 
 \begin{table}[!h]
	\caption{Convergence performance and relative numerical errors by PDDM with different $\delta_S$ and $\delta_D$.}
	\label{OrderDDM}\tabcolsep 0pt \vspace*{-10pt}
	\par
	\begin{center}
		\renewcommand\arraystretch{1.1}
		\def\temptablewidth{1.0\textwidth}
		{\rule{\temptablewidth}{1pt}}
		\begin{tabular*}{\temptablewidth}{@{\extracolsep{\fill}}cc|ccccccccc}
			\hline
			$\delta_D$&$\delta_S$&$h $& $n$ & $\frac{||\mathbf u_S-\mathbf u_{S,h}||_{0}}{||\mathbf u_S||_{0}}$ & $\frac{||\mathbf u_S-\mathbf u_{S,h}||_{1}}{||\mathbf u_S||_{1}}$ & $\frac{||p_S-p_{S,h} ||_{0}}{||p_S||_0}$   & $\frac{||\mathbf u_D-\mathbf u_{D,h}||_{0}}{||\mathbf u_D||_0} $  &  $\frac{||\varphi_D-\varphi_{D,h}||_{0}}{||\varphi_D||_0}$  & $\frac{||\mathbf u_M-\mathbf u_{M,h}||_{0}}{||\mathbf u_M||_0} $  &  $\frac{||\varphi_M-\varphi_{M,h}||_{0}}{||\varphi_M||_0}$\\
			\hline
			1 & 1 & $\frac{1}{8} $&   88    &     0.007027    &     0.055167    &    0.912037    &    0.020925       &   1.107290   &     0.006293      &    9.117670
			\\ 
			& & $\frac{1}{16}$  &134    &      0.001810    &     0.027379   &     0.259521   &     0.005339       &    0.334569    &    0.001663       &    2.64510
			\\ 
			& & $\frac{1}{32}$  &202   &       0.000457 &       0.013659    &    0.074158   &     0.001346      &     0.107216   &     0.000425      &     0.765935
			\\ 
			& & $\frac{1}{64}$  & 301    &      0.000114     &    0.006825   &     0.023994   &     0.000337     &    0.038057   &     0.000107      &    0.251942 
			\\ 
			& & $\frac{1}{128}$  & 444    &      0.000029    &     0.003412   &     0.007124   &     0.000084     &     0.016521   &     0.000027      &     0.075551
			\\ 
			\hline
			1 & $\frac{1}{3}$ & $\frac{1}{8} $&24   &      0.007009    &     0.055160     &   0.207877   &     0.020929     &      0.299250   &     0.006294    &     1.510130
			\\
			& & $\frac{1}{16}$  &20    &      0.001805   &      0.027378    &    0.039257   &     0.005340       &     0.114909   &     0.001663      &     0.235920
			\\
			& & $\frac{1}{32}$  &18  &        0.000455    &     0.013659   &     0.008837   &     0.001346      &     0.054961    &    0.000425     &     0.046765
			\\
			& & $\frac{1}{64}$  &18    &      0.000114   &      0.006825   &     0.002375   &     0.000338       &   0.027271    &    0.000107    &     0.015599
			\\
			& & $\frac{1}{128}$  &18      &    0.000029    &     0.003412   &     0.000698   &     0.000085      &    0.013608    &    0.000027   &     0.007104
			\\
			\hline
		\end{tabular*}%
	\end{center}
\end{table}

\subsection{3D shallow water system with a complicated dual-permeability region}
To present the complicated flow characteristics and to demonstrate the conservation of mass, a 3D model is simulated, then the flow speed, streamlines, and the numerical interface error are depicted. The 3D shallow water system with a complicated dual-permeability region is inspired by \cite{Sun2021, Sun2021tg}, and the details of the domain structure are omitted and can be referred to the above references.

In the shallow water channel, $\mathbf{u}_S = [4y(1-y), 0.0]$ is imposed as the inflow surface velocity and the free boundary condition is defined on the outlet surface. Then we set no-slip boundary conditions for the rest of the boundaries except the interface. As for the complicated dual-permeability region, we assume  no-flux conditions
on the vertical faces and the two impermeable solids faces, and a free boundary condition on the bottom surface. Moreover, let $\nu=1.0,\ \sigma=1.0, \ \mu=1.0, \ \rho=1.0, \ \alpha=1.0, \ k_D=1.0$ and $k_M=0.01$. At the same time, the external body force $\mathbf{f}_S$ and source term $f_D$ are both zero for convenient. We choose $\delta_S = \frac{1}{3}, \delta_D = 1.0$ for rapid computation.

We utilize the proposed PDDM algorithm to compute the above 3D model with mesh size $h=\frac{1}{8}$ and $h=\frac{1}{16}$. In Fig. 5.2, the velocity distribution and some representative streamlines from Stokes channel to the microfracture are shown to illustrate the effectiveness of the proposed algorithm. It can be seen from Fig. 5.2 that the flow speed and direction are consistent with the actual situation and are comparable with the results of \cite{Sun2021, Sun2021tg}. We also present the matrix flow speed and streamlines in Fig. 5.3, which shows that there is no-fluid communication between the matrix and Stokes channel and the matrix velocity is much smaller than microfracture's.
\begin{figure}[htbp]
	\label{SDP3Dwhole_ex3}
	\centering
	\includegraphics[width=50mm,height=60mm]{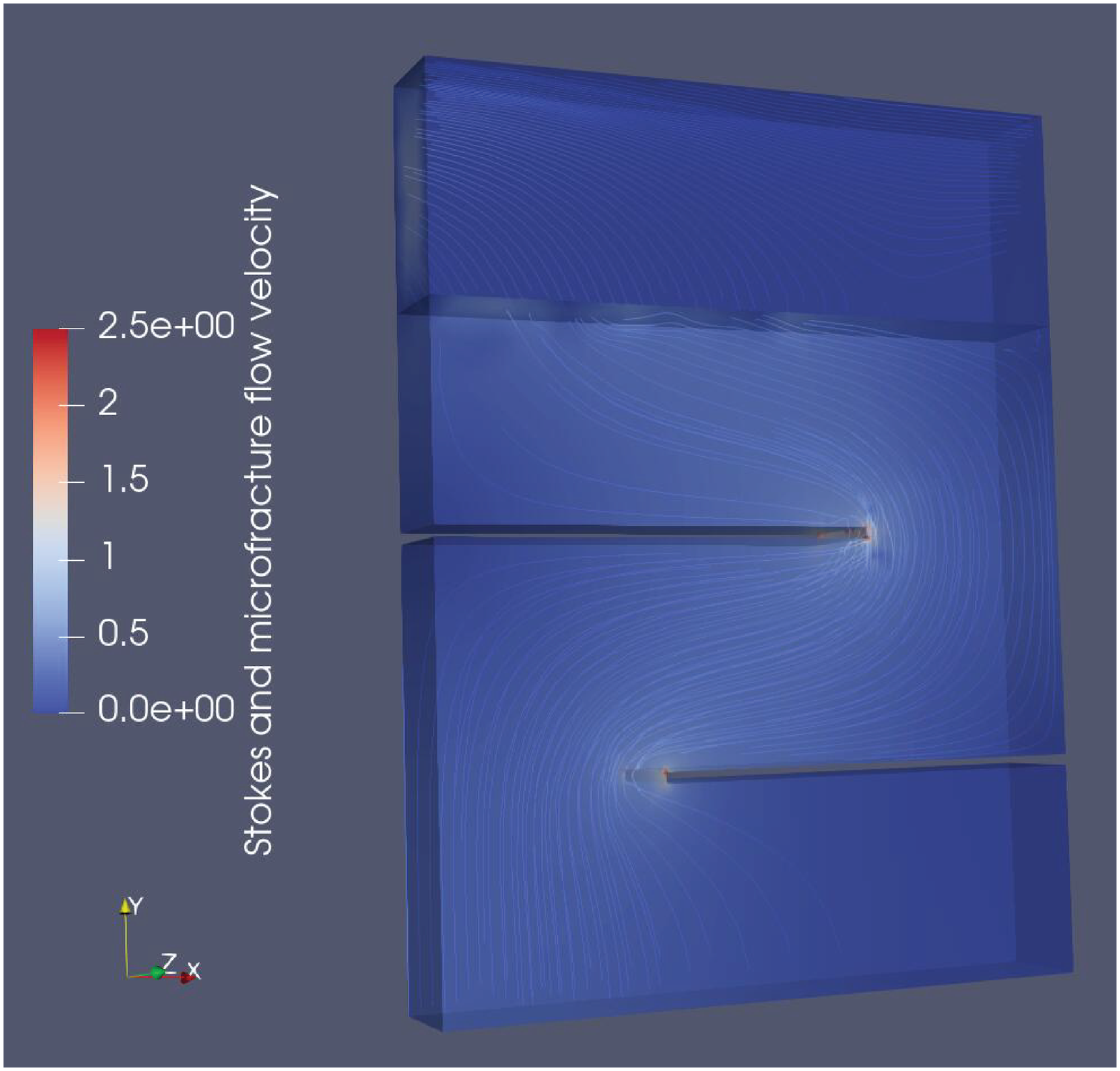}
	\hspace{4mm}
	\includegraphics[width=50mm,height=60mm]{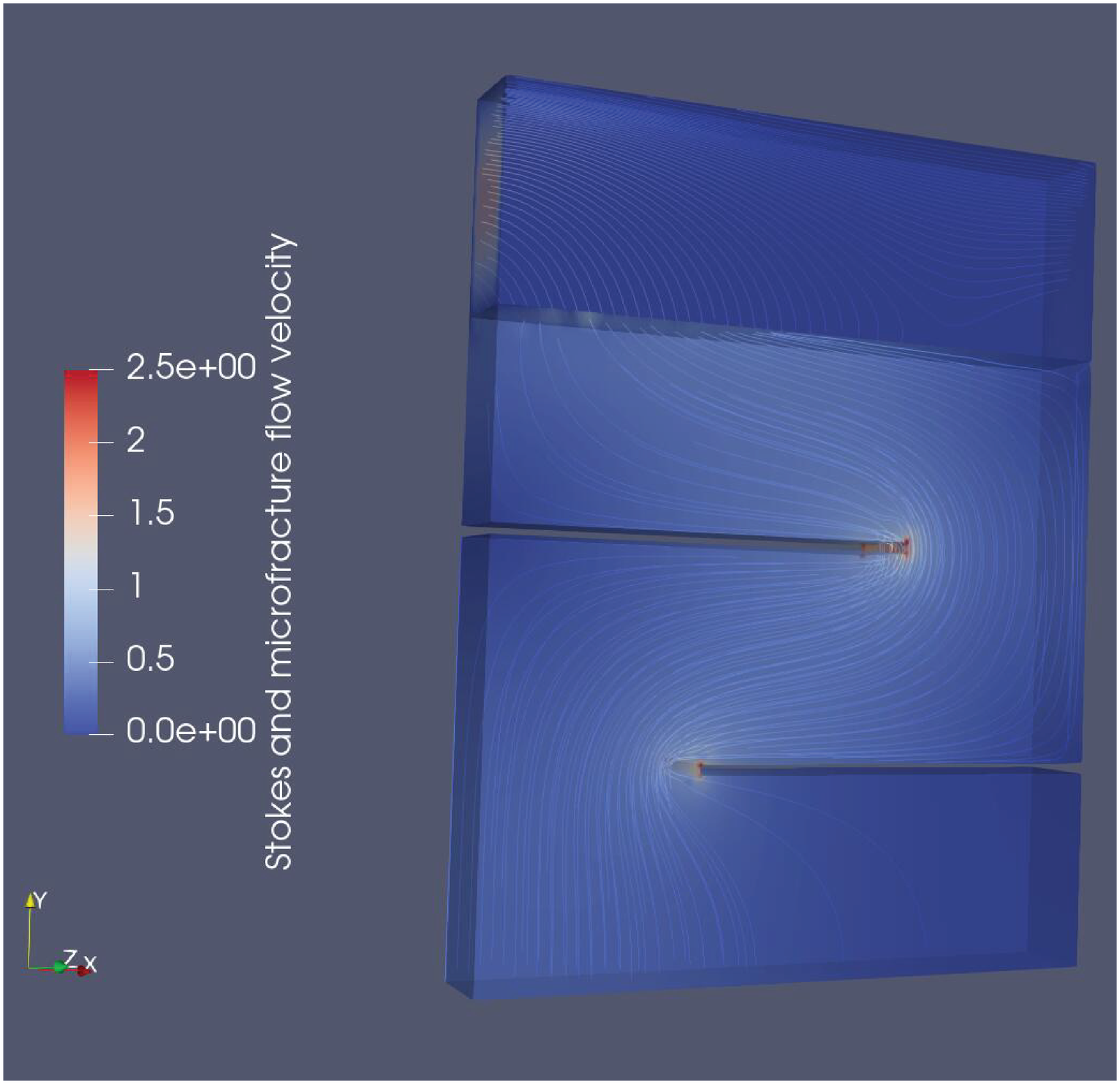}
	\caption{The velocity distribution and some representative streamlines from Stokes channel to the microfracture with mesh size $h=\frac{1}{8}$ (left) and $\frac{1}{16}$ (right).}
\end{figure}
\begin{figure}[htbp]
	\label{SDP3Dmatrix_ex3}
	\centering
	\includegraphics[width=50mm,height=60mm]{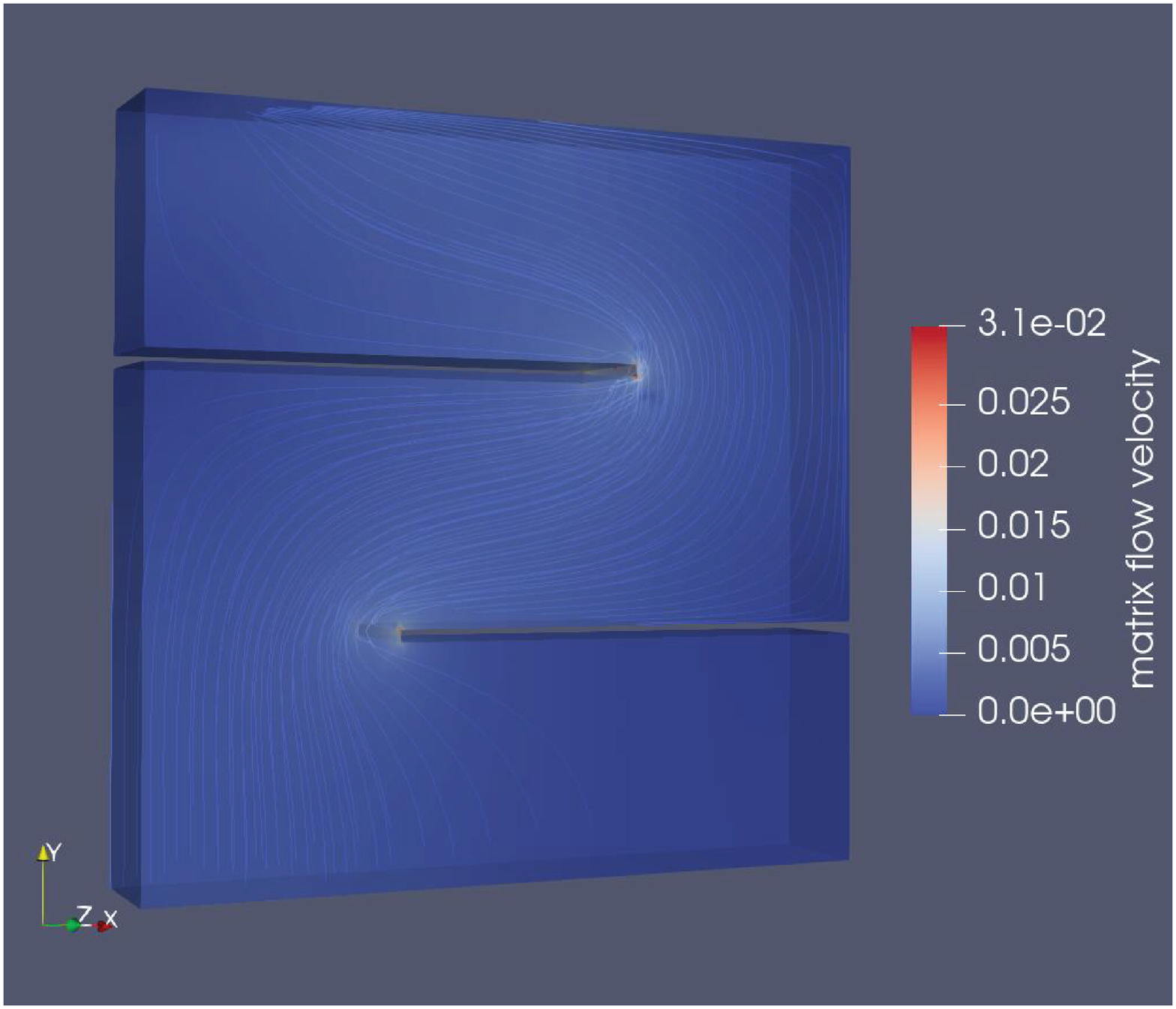}
	\hspace{4mm}
	\includegraphics[width=50mm,height=60mm]{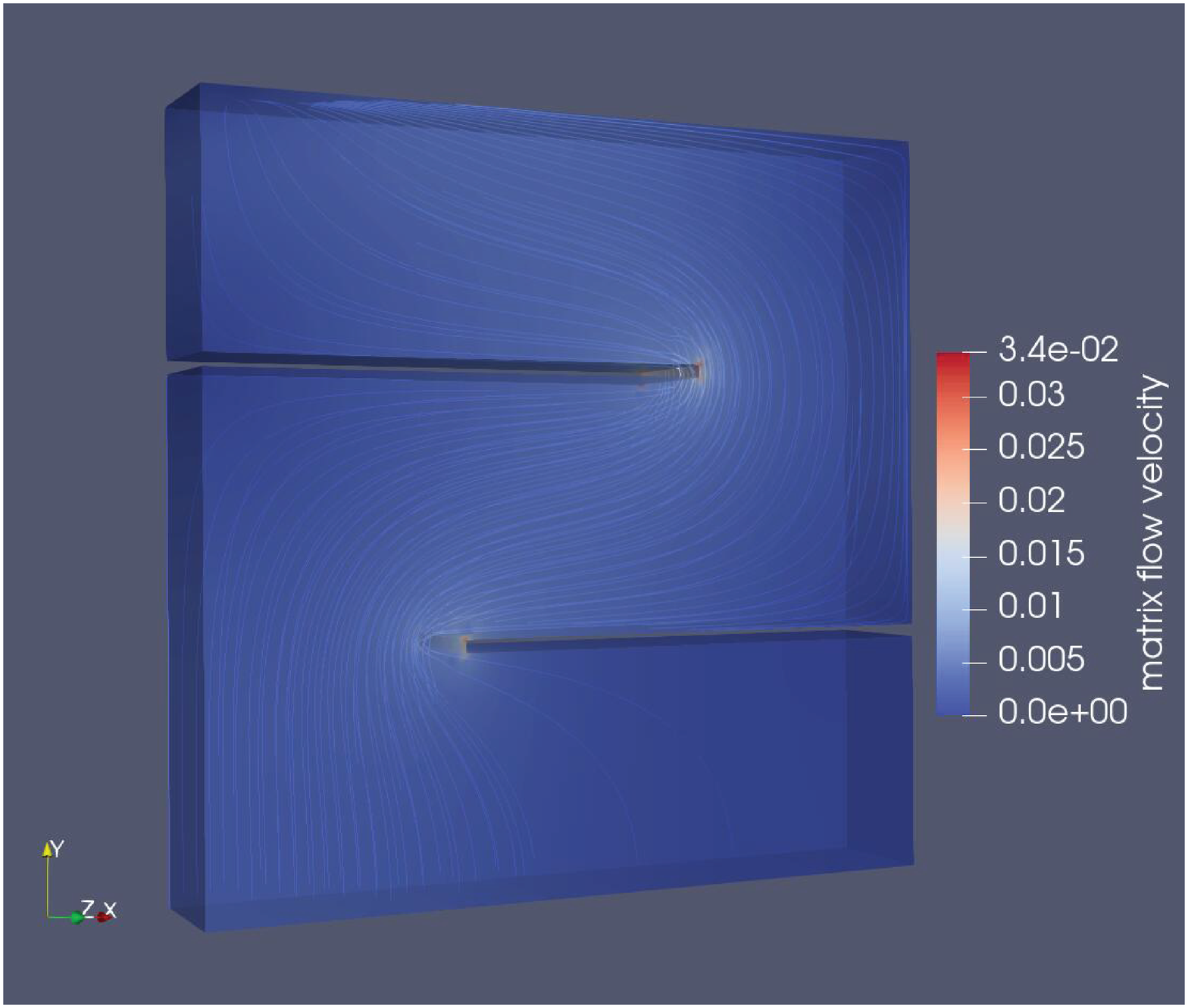}
	\caption{The matrix flow speed and streamlines with mesh size $h=\frac{1}{8}$ (left) and $\frac{1}{16}$ (right).}
\end{figure}

For multi-domain, multi-physics problems, it is necessary to check whether the mass conservation is satisfied on the interface. Hence, we compute the numerical interface errors $|\mathbf{u}_S\cdot\mathbf{n}_S+\mathbf{u}_D\cdot\mathbf{n}_D|$ on the interface cross-section to check whether numerical solutions satisfy the conservation of mass interface condition (\ref{interface1}). We plot the numerical interface mass errors with different mesh scales $h$ in Fig. 5.4, which demonstrates the continuity of the normal velocity on the interface. Furthermore, while the mesh become finer, the errors become smaller as expected.
\begin{figure}[htbp]
	\label{SDP3interface_ex3}
	\centering
	\includegraphics[width=70mm,height=25mm]{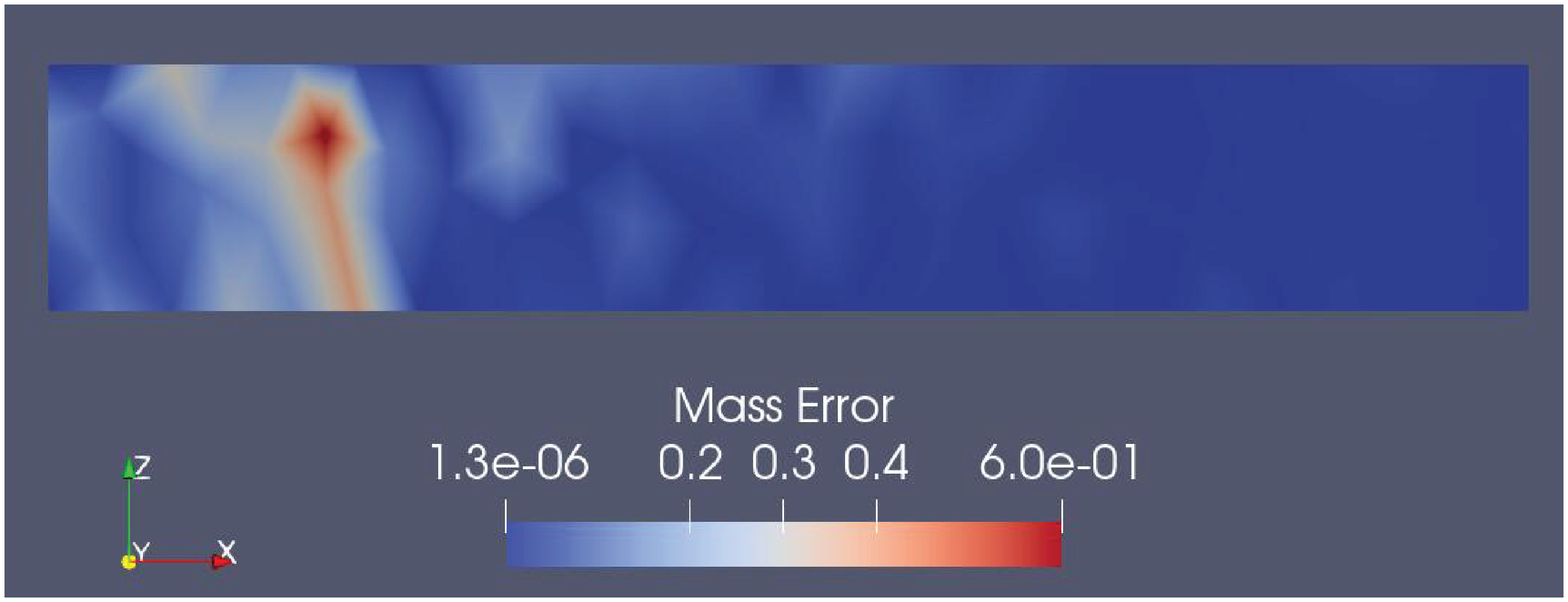}
	\hspace{1mm}
	\includegraphics[width=70mm,height=25mm]{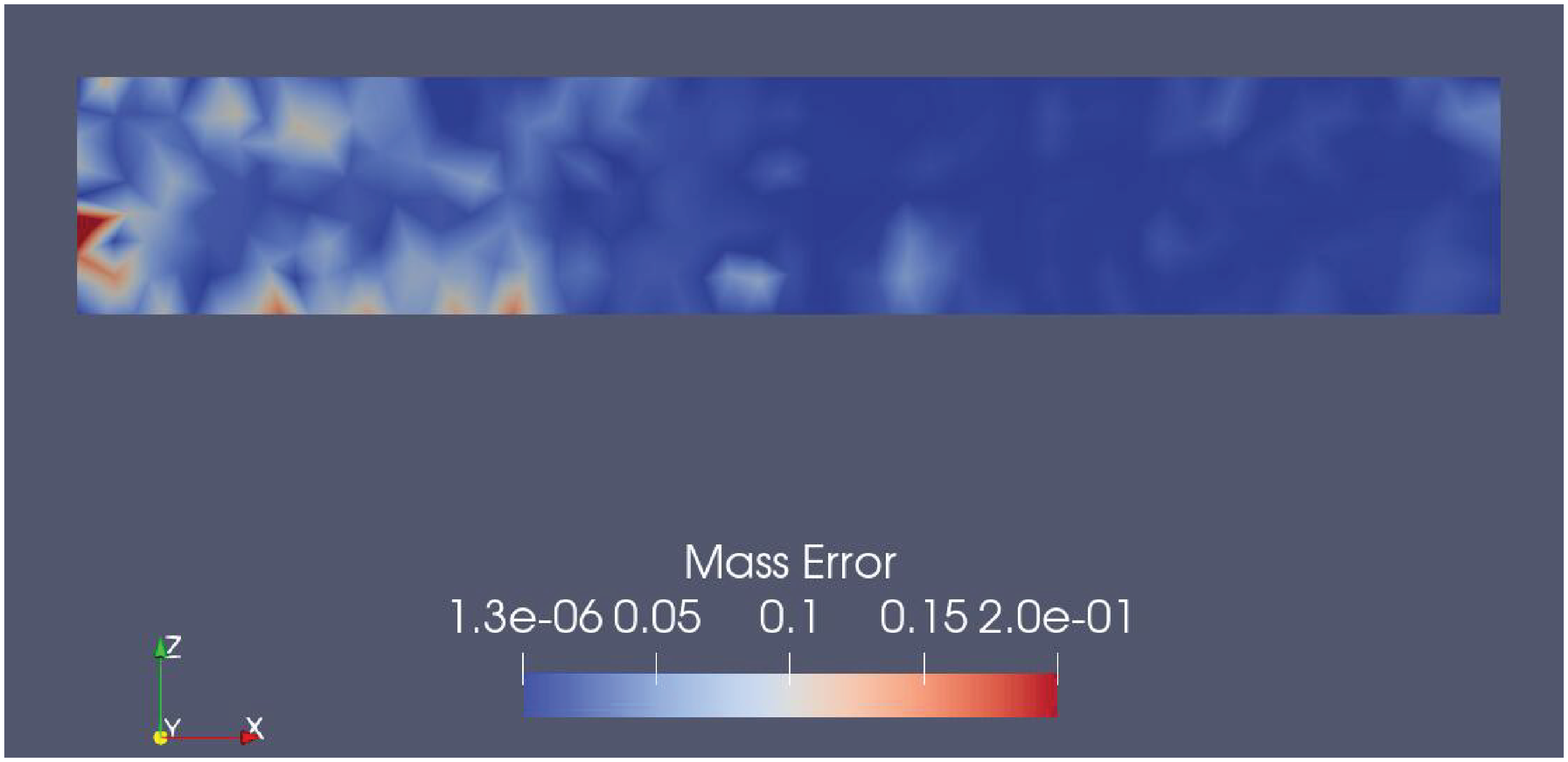}
	\caption{The numerical interface mass errors with mesh size $h=\frac{1}{8}$ (left) and $\frac{1}{16}$ (right).}
\end{figure}

\subsection{Horizontal cased-hole completion wellbore with a vertical production wellbore}
In order to extract oil or gas more efficiently from the unconventional naturally fractured reservoir, one of the most important techniques in petroleum engineering is utilizing the horizontal cased-hole completion wellbore with a vertical production wellbore \cite{Mahbub2019, Mahbub2020}. The main purpose of this experiment is to show the patterns of the flow around the horizontal cased-hole completion wellbore.

We construct a simple geometrical shape to perform this simulation and explain the treatment of boundary conditions, as shown in Fig. 5.5. Eight hydraulic fractures are embedded in the horizontal wellbore as the fluid transmission interface boundaries $\Gamma$. Other interfaces are no communication boundaries $\Gamma_{no}^{SD}$, which is determined by the fundamental properties of the multistage hydraulic fractured horizontal cased-hole completion wellbore, i.e.
\begin{eqnarray*}
\mathbf{u}_S = [0.0,0.0], \hspace{4mm} \mathbf{u}_D \cdot \mathbf{n}_D = 0.0, \hspace{4mm} \mathbf{u}_M \cdot \mathbf{n}_D = 0.0, \hspace{4mm} \mathtt{on} \ \Gamma_{no}^{SD}.
\end{eqnarray*}
The inlet boundary conditions in the dual-permeability domain are imposed as follows:
\begin{eqnarray*}
 \mathbf{u}_D \cdot \mathbf{n}_D = 2.0, \hspace{4mm} \mathbf{u}_M \cdot \mathbf{n}_D = 0.01, \hspace{4mm} \mathtt{on} \ \Gamma_{D}.
\end{eqnarray*}
For the fluid region, we assume a free boundary condition on the outlet boundary $\Gamma_{out}$ and no-slip boundary conditions on the boundary $\Gamma_S$. Moreover, the parameters of this model are selected as $\nu=1.0,\ \sigma=1.0, \ \mu=1.0, \ \rho=1.0, \ \alpha=1.0, \ \mathbf{f}_S = [0.0, 0.0], \ f_D =0.0, \ k_D=10^{-4}$ and $k_M=10^{-8}$.
\begin{figure}[htbp] 
	\label{domain_ex2}
	\centering
	\includegraphics[width=65mm,height=50mm]{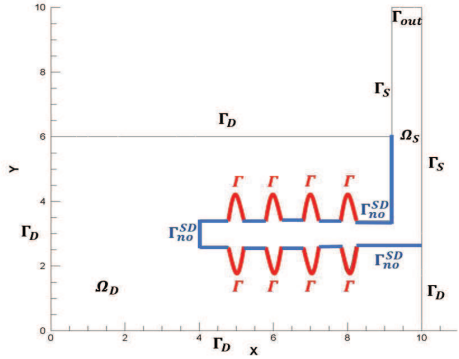}
	\caption{A simple geometrical shape to illustrate the horizontal cased-hole completion wellbore with a vertical production wellbore and the treatment of boundary conditions.}
\end{figure}

We solve this petroleum productivity model by PDDM on a uniform mesh with $h=\frac{1}{10}$. Through some tests, PDDM seems difficult to converge for $\delta_S \leq \delta_D$ while the intrinsic permeabilities $k_D$ and $k_M$ are small. Inspired by the numerical experiments in \cite{Chen, Cao11} and our trying \cite{Shi22}, we find that the Robin parameters $\delta_S>\delta_D$ might provide decent convergence results. So let $\delta_S=5.0$ and $\delta_D=1.0$ in this simulation to get convergence result. We present the pressure field of the microfracture, the multistage hydraulic fractured horizontal wellbore and the matrix in Fig. 5.6. The warmer color indicates higher pressure, so that the higher pressure in the matrix and microfractures can push the fluid into the wellbore with relatively lower pressure. Then, the fluid flow speed and streamlines around the multistage hydraulic fractured horizontal wellbore with a vertical production wellbore are shown in Fig. 5.7. The direction of the flow streamlines supports the expected results, and the deep blue color around the horizontal wellbore means the fluid flow does not interact directly with the $\Gamma_{no}^{SD}$. The regular velocity streamlines across the interface and the smooth distribution of the numerical pressure demonstrate the stability of the proposed DDM.
\begin{figure}[htbp]
	\label{SDPpressure_ex2}
	\centering
	\includegraphics[width=65mm,height=50mm]{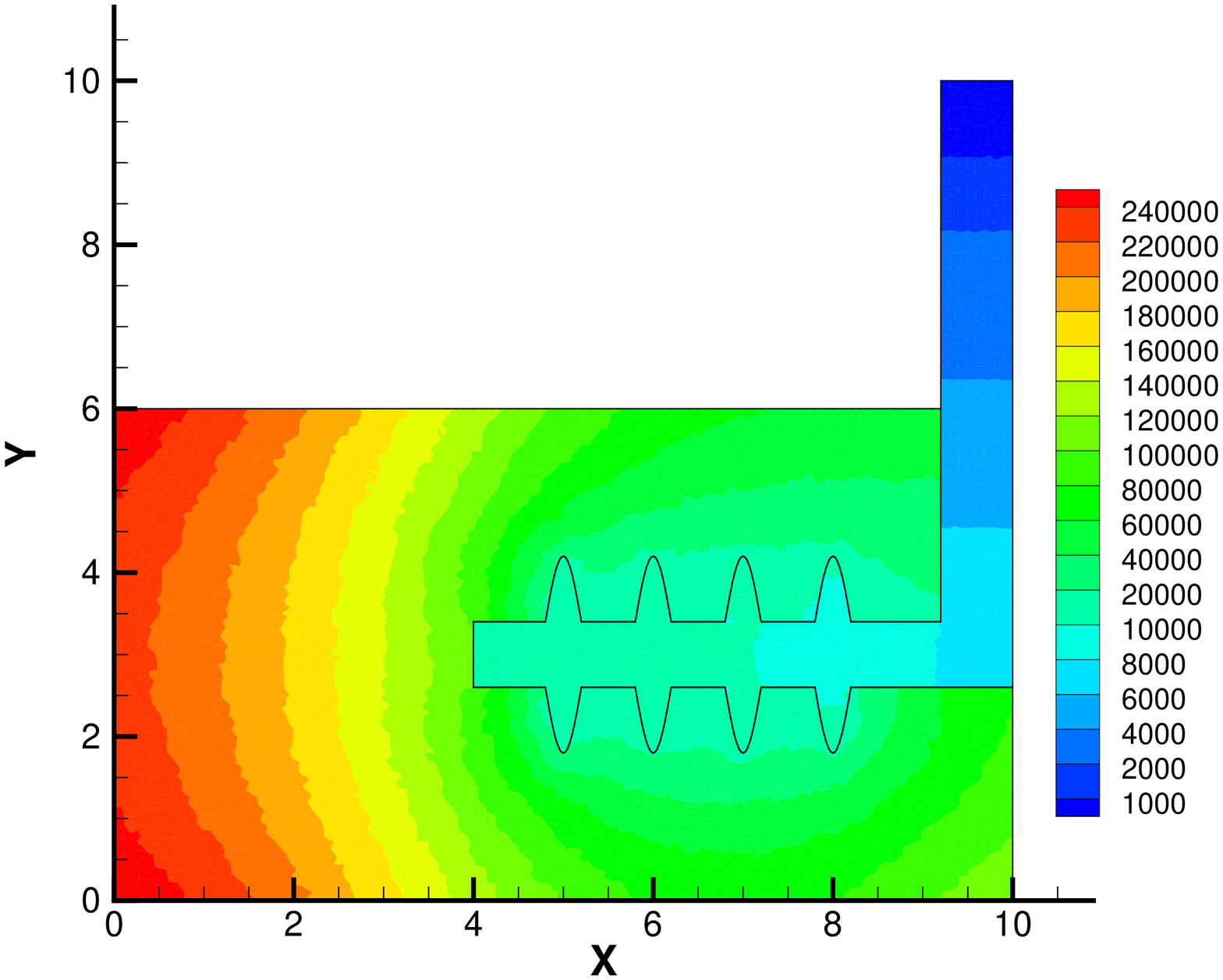}
	\hspace{4mm}
	\includegraphics[width=65mm,height=50mm]{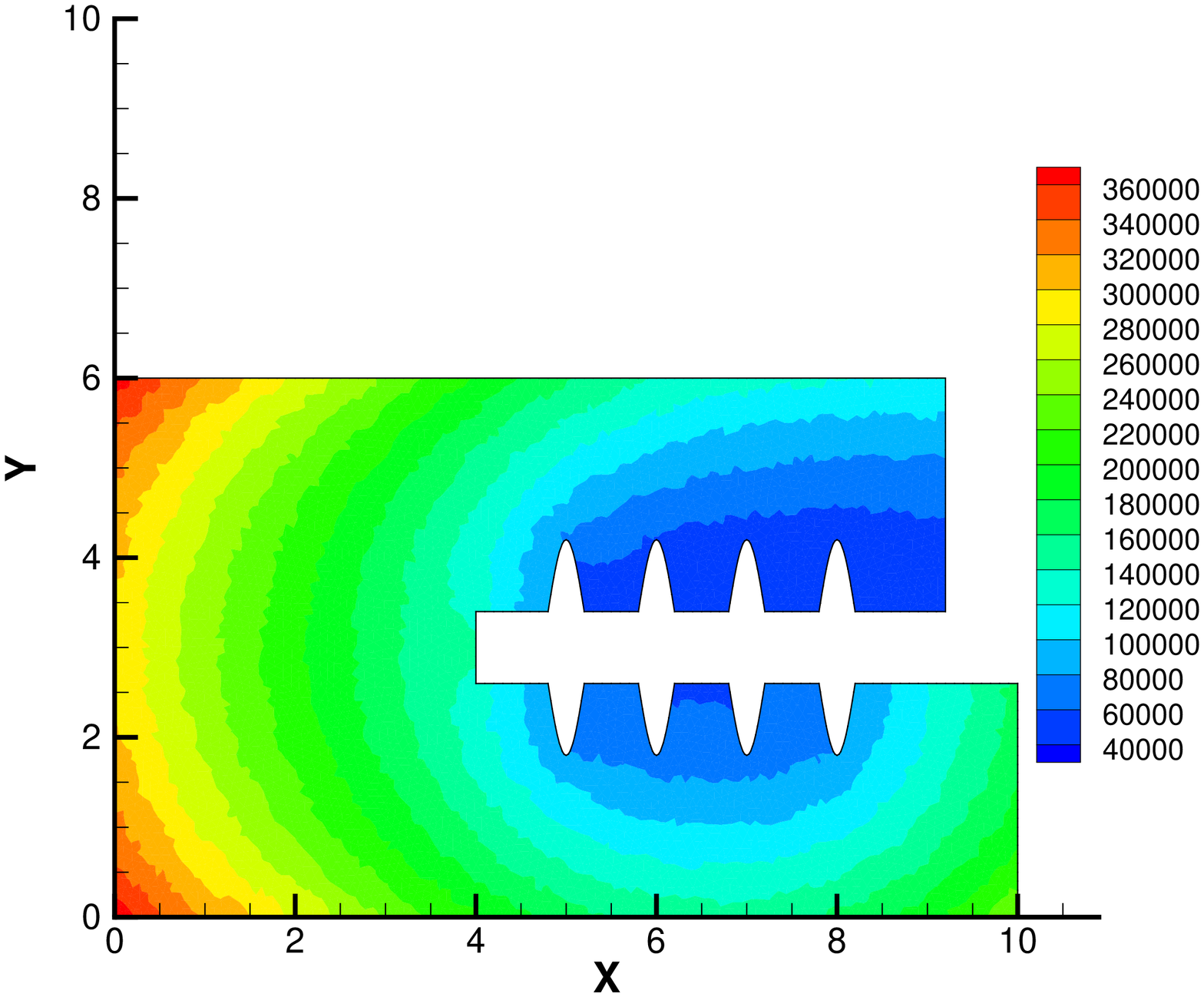}
	\caption{The pressure field of the microfracture, the multistage hydraulic fractured horizontal wellbore (left) and the pressure in the matrix (right).}
\end{figure}
\begin{figure}[htbp]
	\label{SDPvelocity_ex2}
	\centering
	\includegraphics[width=65mm,height=60mm]{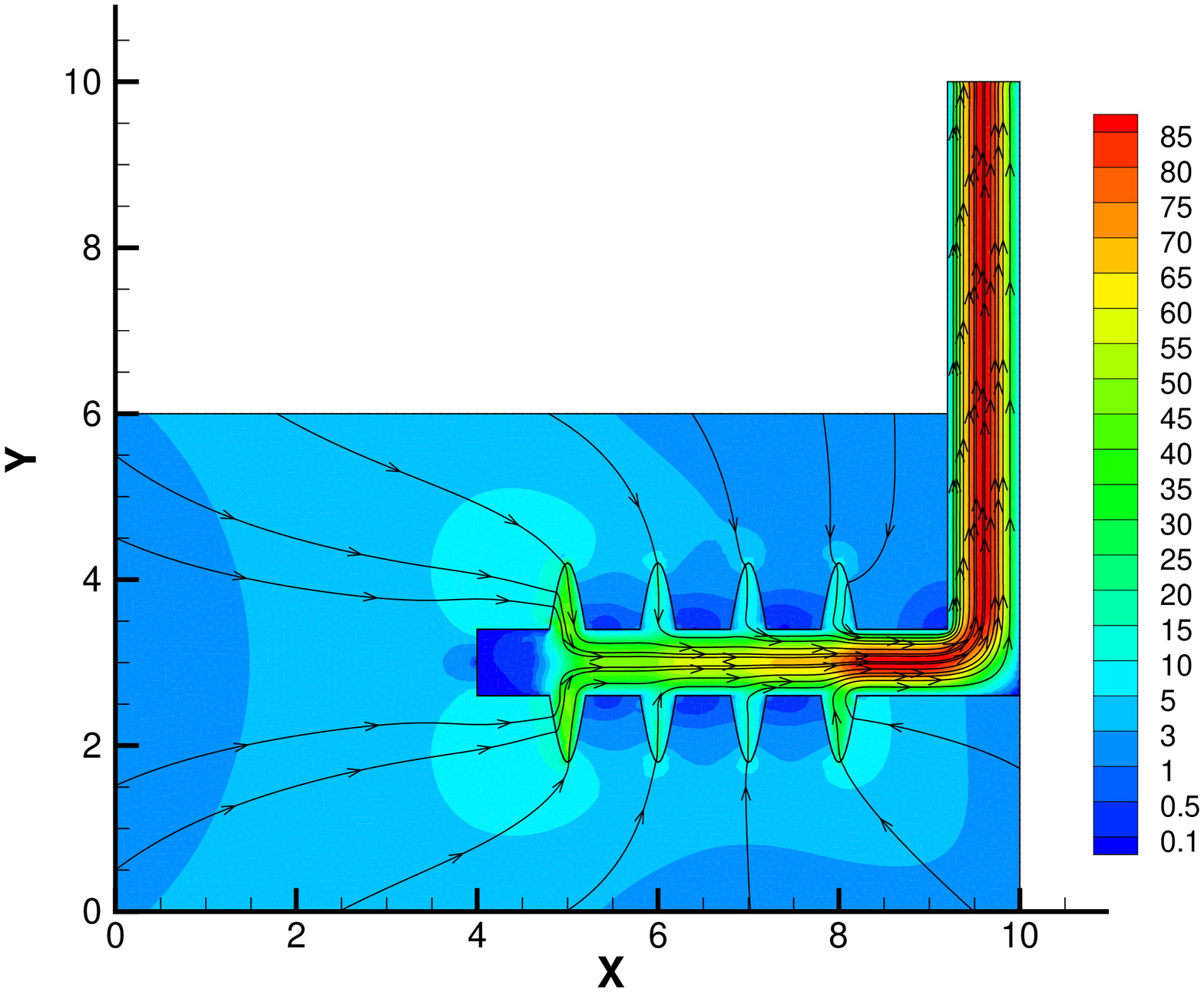}
	\hspace{4mm}
	\includegraphics[width=65mm,height=60mm]{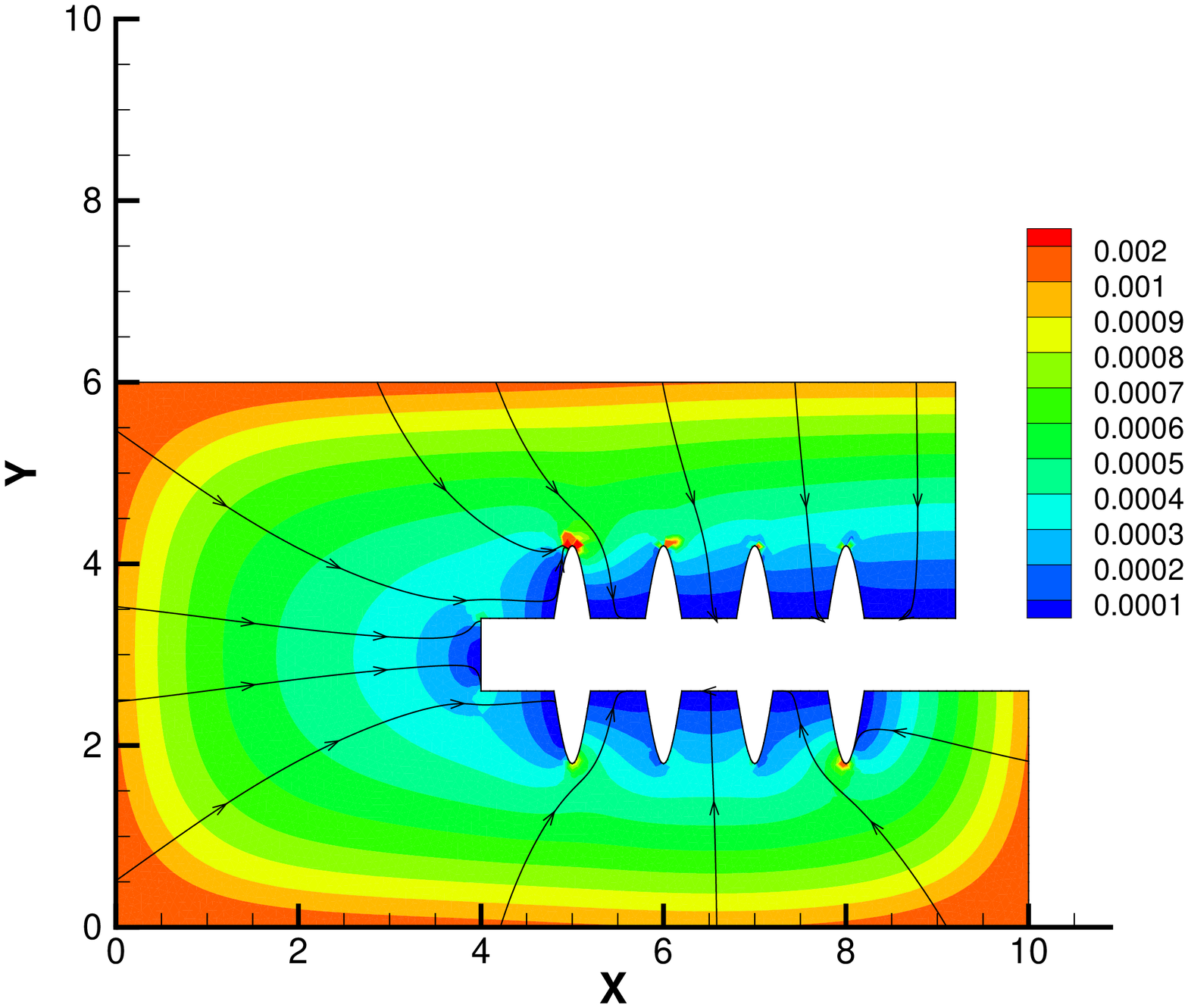}
	\caption{The fluid flow speed and streamlines in the microfracture, the multistage hydraulic fractured horizontal wellbore (left) and in the matrix (right).}
\end{figure}

\section{Conclusions}
In this paper, a novel parallel domain decomposition   (PDDM) algorithm is constructed for solving the fully-mixed Stokes-dual-permeability model with the physically realistic BJ interface conditions. This PDDM can completely decouple the Stokes-dual-permeability model into three independent subproblems, which could be solved in parallel. In order to obtain a strict convergence analysis of PDDM, we propose a general and novel convergence Lemma 4.1. Moreover, by choosing two suitable parameters, $\delta_S<\delta_D$, the convergence rates are independent of the mesh size $h$. The final aim is to extend such PDDM algorithm to the fully-mixed Stokes–dual-permeability system coupled with transport equations, which attaches great importance to the velocities of the microfracture and matrix flow.

\end{document}